\input ./style/arxiv-general.cfg
\documentclass[MSNbibl,number,citesort,dvips]{arxbj}
\makeatletter
   \@ifpackageloaded{graphicx}{}{\usepackage{graphicx}}
\makeatother

\usepackage{accents}

\volume{23}
\issue{1}
\pubyear{2017}
\firstpage{710}
\lastpage{742}
\doi{10.3150/15-BEJ697}
\docsubty{FLA}

\makeatletter
\newcommand{\aaa}{\accentset{\circ}{I}}
\newcommand{\rrvert}{\vert}
\newcommand{\rrVert}{\Vert}
\newcommand{\llvert}{\vert}
\newcommand{\llVert}{\Vert}
\newtheorem{Thm}{Theorem}[section]
\newtheorem{Lem}[Thm]{Lemma}
\newtheorem{Pro}[Thm]{Proposition}
\newproclaim{Def}[Thm]{Definition}
\newremark{Eg}[Thm]{Example}
\newremark{Rem}[Thm]{Remark}
\makeatother

\begin{document}
\begin{frontmatter}

\title{The impact of the diagonals of polynomial forms on limit
theorems with long memory}
\runtitle{Impact of diagonals of polynomial forms on limit theorems
with long memory}

\begin{aug}
\author{\inits{S.}\fnms{Shuyang}~\snm{Bai}\thanksref{e1}\ead[label=e1,mark]{bsy9142@bu.edu}}
\and
\author{\inits{M.S.}\fnms{Murad S.}~\snm{Taqqu}\corref{}\thanksref{e2}\ead[label=e2,mark]{bumastat@gmail.com}}
\address[]{Department of Mathematics and Statistics, Boston
University, 111 Cumminton Street,
Boston, MA 02215, USA. \printead{e1}; \printead*{e2}}
\end{aug}

%
\received{\smonth{3} \syear{2014}}
%
\revised{\smonth{11} \syear{2014}}

%
\begin{abstract}
We start with an i.i.d.
sequence and consider the product of two polynomial-forms moving
averages based on that sequence. The coefficients of the polynomial
forms are asymptotically slowly decaying homogeneous functions so that
these processes have long memory. The product of these two polynomial
forms is a stationary
nonlinear process. Our goal is to obtain limit theorems for the
normalized sums of this product process in three cases: exclusion of
the diagonal terms of the polynomial form,
inclusion, or the mixed case (one polynomial form excludes the
diagonals while the other one includes them). In any one of these
cases, if the product has long memory,
then the limits are given by Wiener chaos. But the limits in each of
the cases are quite different. If the diagonals are excluded, then the
limit is expressed as in the product
formula of two Wiener--It\^o integrals. When the diagonals are
included, the limit stochastic integrals are typically due to a single
factor of the product, namely the one with
the strongest memory. In the mixed case, the limit stochastic integral
is due to the polynomial form without the diagonals irrespective of the
strength of the memory.
\end{abstract}

%
\end{frontmatter}

\section{Introduction}

Let $X(n)$ be a stationary process with mean $0$ and finite variance.
We are interested in the following weak convergence of normalized
partial sum to a process $Z(t)$:
%
\begin{equation}
\label{eqgenerallimit}
\frac{1}{A(N)}\sum_{n=1}^{[Nt]}X(n)
\Rightarrow Z(t)
\end{equation}
as $N\rightarrow\infty$ where $A(N)\rightarrow\infty$ is a suitable
normalization. The limit $Z(t), t\ge0$ if it exists, has stationary
increments and is self-similar with some index $H>0$, that is, for any
$a>0$, $\{Z(at),t\ge0\}$ and $\{a^HZ(t),t\ge0\}$ have the same
finite-dimensional distributions. The parameter $H$ is called the \emph
{memory parameter}\footnote{A precise definition of memory parameter is
given in Definition~\ref{DefHurstexpXn}.} of the process $X(n)$
and the \emph{Hurst index} or \emph{self-similarity parameter} of the
limit process $Z(t)$. The higher the value of $H$, the stronger the
memory of the process $X(n)$.

When the dependence in $X(n)$ is weak, one typically ends up in (\ref
{eqgenerallimit}) with
\[
A(N)= \Biggl(\operatorname{Var}\Biggl[\sum_{n=1}^N
X(n)\Biggr] \Biggr)^{1/2}\sim cN^{1/2}
\]
as $N\rightarrow\infty$ for some $c>0$, and $Z(t)$ is the Brownian
motion. These types of limit theorems are often called \emph{central
limit theorems}.

When, however, the dependence in $X(n)$ is so strong that
$\operatorname{Var}[\sum_{n=1}^N X(n)]$ grows faster than the linear speed $N$, and typically
as $N^{2H}$ with $H\in(1/2,1)$, the limit process $Z(t)$ in (\ref
{eqgenerallimit}) is no longer Brownian motion. $Z(t)$ is in this
case a self-similar process with stationary increments which has a
\emph
{Hurst index} $H$ (see \cite{embrechtsmaejima2002selfsimilar}).
This type of limit theorems involving non-Brownian limits are often
called \emph{noncentral limit theorems}. When the process $X(n)$ is
nonlinear and has long memory, the limit $Z(t)$ can be non-Gaussian
(e.g., \cite{dobrushinmajor1979non,taqqu1979convergence,surgailis1982zones}).

In \cite{baitaqqu2013generalized}, a noncentral limit theorem is
established for an off-diagonal polynomial-form process called $k$th
order \emph{discrete chaos process}:
%
\begin{equation}
\label{eqXnintro1}
Y'(n)=\sum_{0<i_1,\ldots,i_k<\infty}'
a(i_1,\ldots,i_k)\varepsilon _{n-i_1}\cdots
\varepsilon_{n-i_k},
\end{equation}
where the prime $'$ indicates that we do not sum on the diagonals
$i_p=i_q$, $p\neq q$, the noise $\varepsilon_i$'s are i.i.d. random
variables with mean $0$ and variance $1$, and $a(\cdot)$ is
asymptotically some homogeneous function $g$ called \emph{generalized
Hermite kernel} (GHK). The limit $Z(t)$, called a \emph{generalized
Hermite process}, is expressed by a $k$-fold \emph{Wiener--It\^o integral}:
%
\begin{equation}
\label{eqgenhermproc} Z(t)=\int_{\mathbb{R}^k}'\!\int
_0^t g(s-x_1,\ldots,s-x_k)1_{\{
s>x_1,\ldots,s>x_k\}}
\,\mathrm{d}s B(\mathrm{d}x_1)\cdots B(\mathrm{d}x_k),
\end{equation}
where the prime $'$ indicates that we do not integrate on the diagonals
$x_p=x_q$, $p\neq q$, and $B(\cdot)$ is Brownian motion. These
processes $Z(t)$ include the \emph{Hermite process} considered in
\cite{dobrushinmajor1979non,taqqu1979convergence} and \cite{surgailis1982zones}.

In \cite{baitaqqu2014convergence}, a noncentral limit theorem is
established for a polynomial-form process called $k$th order \emph
{discrete Volterra process}:
%
\begin{equation}
\label{eqXnintro} Y(n)=\sum_{0<i_1,\ldots,i_k<\infty} a(i_1,
\ldots,i_k)\varepsilon _{n-i_1}\cdots\varepsilon_{n-i_k},
\end{equation}
which differs from $Y'(n)$ in (\ref{eqXnintro1}) by including the
diagonals, and where $a(\cdot)$ is asymptotically $g(\cdot)$, some
special type of generalized Hermite kernel called \emph{generalized
Hermite kernel of Class (B)} (GHK(B)). The limit $Z(t)$ can be
heuristically thought as (\ref{eqgenhermproc}) with diagonals
included, and is precisely expressed as a $k$-fold \emph{centered
Wiener--Stratonovich integral}, which is a linear combination of
certain Wiener--It\^o integrals of orders lower than or equal to $k$
(see \cite{baitaqqu2014convergence}).

In this paper, we contrast the effect of two types of stationary
sequences in the limit theorem~(\ref{eqgenerallimit}). The first
stationary sequence is
%
\begin{equation}
\label{eqprodofchaos}
X(n)=Y'_1(n)Y'_2(n),
\end{equation}
that is, a product of two long memory chaos processes (\ref{eqXnintro1}) which \emph{exclude} the diagonals. The second stationary
sequence is
%
\begin{equation}
\label{eqprodofvolt} X(n)=Y_1(n)Y_2(n),
\end{equation}
that is, a product of two long memory processes in (\ref{eqXnintro}) which \emph{include} the diagonals. We also consider the mixed case
%
\begin{equation}
\label{eqprodofmixed} X(n)=Y_1'(n)Y_2(n).
\end{equation}

Limit theorems for such types of product are of interest, for example,
in statistical inference involving long memory processes with different
memory parameters (\cite{koulbaillie2004regression}, see also
Proposition~11.5.6 of \cite{giraitiskoulsurgailis2009large}), and
in the study of covariation of fractional Brownian motions with
different Hurst indexes \cite{maejimatudor2012selfsimilar}.
Typically, the factor processes $Y$ there are assumed to be either
linear (or Gaussian) or a transformation of linear process (or
Gaussian), which yields in the limit a generalized Rosenblatt processes
where $g(x_1,x_2)=x_1^{\gamma_1}x_2^{\gamma_2}$ in (\ref{eqgenhermproc}). By taking the factors $Y$ to be some nonlinear processes as in
(\ref{eqprodofchaos}), (\ref{eqprodofvolt}) and (\ref{eqprodofmixed}), one can obtain much richer limit structures, which are briefly
described below.

We show that in the case (\ref{eqprodofchaos}), the limit in (\ref
{eqgenerallimit}) is expressed as Wiener--It\^o integrals which can
be obtained by using a rule similar to that used for computing the
product of two Wiener--It\^o integrals. In fact, if the stationary
sequences $Y'_1(n)$ and $Y'_2(n)$ have, respectively, memory parameters
$H_1,H_2\in(1/2,1)$ with $H_1+H_2>3/2$, then the limit in (\ref
{eqgenerallimit}) has Hurst index
\[
H=H_1+H_2-1\in(1/2,1).
\]
In the case (\ref{eqprodofvolt}), in contrast, the limit stochastic
integrals are typically due to a single factor $Y_1(n)$ or $Y_2(n)$,
namely, the one with the strongest memory parameter. The Hurst index of
the limit is then
\[
\max(H_1,H_2)\in(1/2,1)
\]
which is always greater than $H_1+H_2-1$. In the case (\ref{eqprodofmixed}), only the off-diagonal factor $Y_1'(n)$ contributes to the
limit stochastic integral, irrespective of the strength of the memory.

The paper is organized as follows. Section~\ref{secbackground}
contains some background. We state the main results in Section~\ref{secstateresult}, namely, Theorem~\ref{Thmgeneralprodchaos} for
processes without diagonals, Theorem~\ref{Thmgeneralprodvolt} for
processes with diagonals and Theorem~\ref{Thmgeneralprodmixed} for
the mixed case. Section~\ref{secprelim} provides some preliminary
results used in the proofs. Section~\ref{secproof} contains the proofs
of the theorems.

\section{Background}\label{secbackground}
The following notation will be used throughout: $\mathbf{0}$ denotes
the zero vector $(0,0,\ldots,0)$ and $\mathbf{1}=(1,1,\ldots,1)$
denotes the vector with ones in every component. For two vectors
$\mathbf{x}$ and $\mathbf{y}$ with the same dimension, we write
$\mathbf
{x}\le\mathbf{y}$ (or $<$, $\ge$, $>$) if the inequality holds
componentwise. We let
\[
[x]=\sup\{n\in\mathbb{Z}\dvt n\le x\}
\]
for any real $x$ and for a real vector $\mathbf{x}=(x_1,\ldots,x_k)$,
we define
\[
[\mathbf{x}]=\bigl([x_1],\ldots,[x_k]\bigr).
\]
The notation $1_A$ denotes the indicator function of a set $A$. The
value of a constant $C>0$ or $c>0$ may change from line to line.

In \cite{baitaqqu2013generalized}, the following classes of
functions were introduced.

\begin{Def}\label{DefGHK}
A measurable function $g$ defined on $\mathbb{R}_+^k$ is called a
\emph
{generalized Hermite kernel} (GHK) with homogeneity exponent
%
\begin{equation}
\label{eqalpha} \alpha\in \biggl(-\frac{k+1}{2},-\frac{k}{2} \biggr),
\end{equation}
if it satisfies
\begin{enumerate}[2.]
\item[1.] $g(\lambda\mathbf{x})=\lambda^{\alpha}g(\mathbf{x})$,
$\forall
\lambda>0$;
\item[2.] $\int_{\mathbb{R}_+^k} \llvert g(\mathbf{1}+\mathbf{x})g(\mathbf
{x})\rrvert \,\mathrm{d}\mathbf{x}<\infty$;
\end{enumerate}
A GHK $g$ is said to belong to Class (B) [abbreviated as GHK(B)], if
$g$ is a.e. continuous on $\mathbb{R}_+^k$ and
\[
\bigl\llvert g(\mathbf{x})\bigr\rrvert \le c \llVert \mathbf{x}\rrVert
^{\alpha}=c(x_1+\cdots +x_k)^{\alpha}
\]
($\|\cdot\|$ is the $L^1$-norm) for some constant $c>0$.
\end{Def}

\begin{Rem}\label{RemGHKhtL^2}
As it was shown in Theorem~3.5 of \cite{baitaqqu2013generalized},
if $g$ is a GHK, then
\[
\int_0^t \bigl\llvert g(s\mathbf{1}-
\mathbf{x})\bigr\rrvert 1_{\{s\mathbf{1}>\mathbf{x}\}
}\,\mathrm{d}s<\infty
\]
for a.e. $\mathbf{x}\in\mathbb{R}^k$, and the function
\[
h_t(\mathbf{x}):=\int_0^t g(s
\mathbf{1}-\mathbf{x})1_{\{s\mathbf
{1}>\mathbf{x}\}}\,\mathrm{d}s\in L^2\bigl(
\mathbb{R}^k\bigr).
\]
\end{Rem}

Using a GHK, one can define a self-similar process with stationary
increments on a Wiener chaos as follows.

\begin{Def}\label{DefGHprocess}
Let $g$ be a GHK on $\mathbb{R}_+^k$ with homogeneity exponent $\alpha
\in(-\frac{k+1}{2},-\frac{k}{2})$, then (\ref{eqgenhermproc}) is
called a \emph{generalized Hermite process} $Z(t)$. It is self-similar
with Hurst index
%
\begin{equation}
\label{eqH}
H=\alpha+k/2+1.
\end{equation}
\end{Def}

\begin{Eg}\label{EgHermKernel}
If
\[
g(\mathbf{x})=\prod_{j=1}^k
x_j^{\gamma},
\]
where $-1/2-1/k<\gamma<-1/2$, then $Z(t)$ in (\ref{eqgenhermproc})
is the \textit{Hermite process} considered in \cite{dobrushinmajor1979non} and~\cite{taqqu1979convergence}.
\end{Eg}

Note that GHK(B) does not include the kernel in Example~\ref{EgHermKernel}. We use a GHK(B) because of its boundedness property. The
subclass of GHK(B) is, in fact, a dense subset in the whole class of
GHK (see Remark~3.17 of \cite{baitaqqu2013generalized}).

We now state two limit theorems, the first for the discrete chaos
process $Y'(n)$ defined in (\ref{eqXnintro1}) where the diagonals
are excluded, and the second for the Volterra process $Y(n)$ defined in
(\ref{eqXnintro}) which includes the diagonals.

Suppose that $g$ is a GHK(B) on $\mathbb{R}_+^k$, $L(\cdot)$ is a
bounded function defined on $\mathbb{Z}_+^k$ such that
\[
\lim_{n\rightarrow\infty}L\bigl([n\mathbf{x}]+\mathbf{B}(n)\bigr)=1
\]
for any $\mathbf{x}\in\mathbb{R}_+^k$ and any $\mathbb{Z}_+^k$-valued
bounded function $\mathbf{B}(n)$, and suppose that the coefficient
$a(\cdot)$ in (\ref{eqXnintro1}) is given by
%
\begin{equation}
\label{eqa=gL}
a(\mathbf{i})=g(\mathbf{i})L(\mathbf{i}).
\end{equation}

\begin{Pro}[(Theorem~6.5 of \cite{baitaqqu2013generalized})]\label{Thmncltchaossingle}
The following weak convergence holds in $D[0,1]$:
%
\begin{equation}
\frac{1}{N^{H}} \sum_{n=1}^{[Nt]}
Y'(n) \Rightarrow Z(t):=I_k(h_t),
\end{equation}
where $H=\alpha+k/2+1\in(1/2,1)$,
%
\begin{equation}
\label{eqht}
h_t(\mathbf{x})=\int_0^t
g(s\mathbf{1}-\mathbf{x})1_{\{s\mathbf
{1}>\mathbf{x}\}}\,\mathrm{d}s
\end{equation}
with $g$ as in (\ref{eqa=gL}), and $I_k(\cdot)$ denotes the $k$-fold
Wiener--It\^o integral, so that $Z(t)$ is a generalized Hermite process
(\ref{eqgenhermproc}).
\end{Pro}

We now consider the limit when the diagonals are included. If $g$ is
GHK(B) on $\mathbb{R}_+^k$ and is in addition symmetric, we define the
following function $g_r$ by identifying $r$ pairs of variables of $g$
and integrating them out, as follows:
%
\begin{equation}
\label{eqgr} g_r(\mathbf{x})=\int_{\mathbb{R}_+^r}g(y_1,y_1,
\ldots ,y_r,y_r,x_1,\ldots ,x_{k-2r})\,\mathrm{d}
\mathbf{y}.
\end{equation}
In \cite{baitaqqu2014convergence}, a noncentral limit theorem was
established for the Volterra process $Y(n)$ in (\ref{eqXnintro}).
Let
\[
a(\cdot)=g(\cdot)L(\cdot)
\]
in (\ref{eqXnintro}) be given as in (\ref{eqa=gL}) assuming in
addition that $g$ is symmetric.

\begin{Pro}[(Theorem~6.2 of \cite{baitaqqu2014convergence})]\label
{Thmncltvoltsingle}
One has the following weak convergence\vspace*{-4pt} in $D[0,1]$:
%
\begin{equation}
\label{eqncltvolt} \frac{1}{N^{H}} \sum_{n=1}^{[Nt]}
Y(n) \Rightarrow Z(t):=\sum_{0\le
r<k/2} d_{k,r}
Z_{k-2r}(t),
\end{equation}
where\vspace*{-4pt} $H=\alpha+k/2+1\in(1/2,1)$,
%
\begin{equation}
\label{eqdkr} d_{k,r}=\frac{k!}{2^r(k-2r)!r!},
\end{equation}
and\vspace*{-4pt}
%
\begin{equation}
\label{eqZt=gr} Z_{k-2r}(t):=\int'_{\mathbb{R}^{k-2r}}
\!\int_0^t g_r(s\mathbf {1}-\mathbf
{x})\mathrm{1}_{\{s\mathbf{1}>\mathbf{x}\}}\,\mathrm{d}s B(\mathrm{d}x_1)\cdots B(\mathrm{d}x_k)
\end{equation}
is a $(k-2r)$th order generalized Hermite process with GHK given by
$g_r$ in (\ref{eqgr}).
\end{Pro}

\begin{Rem}
The limit process $Z(t)$ in (\ref{eqncltvolt}) can be simply
expressed in terms of a centered Wiener--Stratonovich integral
$\accentset{\circ}{I}_k^c(\cdot)$ as
%
\begin{equation}
\label{eqwienerstratonintegral} Z(t)=\aaa_k^c(h_t),
\end{equation}
where $h_t$ is as in (\ref{eqht}), and\vspace*{-3pt} where
\[
\aaa_k^c(\cdot)=\sum
_{0\le r<k/2} d_{k,r}I_{k-2r}\bigl(
\tau^r \cdot\bigr).
\]
The integral $\aaa_k^c(\cdot)$ differs from the
Wiener--Stratonovich integral
\[
\aaa_k(\cdot):=\sum_{0\le r \le[k/2]}
d_{k,r}I_{k-2r}\bigl(\tau^r \cdot\bigr)
\]
introduced in \cite{humeyer1988integrales} by excluding the term
$r=k/2$ when $k$ is even. Here, the operator $\tau^r$ identifies~$r$
pairs of variables of $h$ and integrates them out (see \cite{baitaqqu2014convergence}). The operator $\tau^r$ is often called a
``trace operator.''\vspace*{-2pt}
\end{Rem}

\section{Statement of the main results}\label{secstateresult}
We state here the main results, and defer the proofs to Sections~\ref{secproofchaos} and \ref{secproofvolt}.
In the statement of the results, the following expressions are used.

\begin{Def}\label{Defcltnclt}
Let $X(n)$ be a stationary process with finite variance. We say that:
\begin{enumerate}[2.]
\item[1.] $X(n)$ satisfies a \emph{central limit theorem} (CLT),\vspace*{-2pt} if
%
\begin{equation}
\label{eqclt} N^{-1/2}\sum_{n=1}^{[Nt]}
\bigl[X(n)-\mathbb{E}X(n)\bigr]\Rightarrow\sigma B(t)
\end{equation}
in $D[0,1]$,
where\vadjust{\goodbreak} $\sigma^2=\sum_{n=-\infty}^\infty\operatorname{Cov}(X(n),X(0))$;

\item[2.] $X(n)$ satisfies a \emph{noncentral limit theorem} (NCLT) with a
Hurst index $H\in(1/2,1)$ and limit $Z(t)$, if
%
\begin{equation}
\label{eqnclt} N^{-H}\sum_{n=1}^{[Nt]}
\bigl[X(n)-\mathbb{E}X(n)\bigr]\Rightarrow Z(t)
\end{equation}
in $D[0,1]$.
\end{enumerate}
\end{Def}

\begin{Rem}
In case~1 above, the ``long-run variance'' $\sigma^2$ can be $0$. In
this case, we understand the limit theorem as degenerate (the
normalization $N^{-1/2}$ is too strong). We do not consider here limit
theorems involving a Hurst index $H<1/2$. In case~2, the limit in (\ref
{eqnclt}) may be fractional Brownian motion.
\end{Rem}

We now consider separately the cases where the diagonals of the
polynomial forms are excluded (chaos processes) and when they are
included (Volterra processes).
\subsection{Limit theorem for a product of long-memory chaos processes}
Suppose that we have the following two discrete chaos processes
(off-diagonal polynomial forms):
%
\begin{equation}
\label{eqoff-diagonalY} Y_1'(n)=\sum
_{\mathbf{i}\in\mathbb{Z}_+^{k_1}}'a^{(1)}(\mathbf {i})
\varepsilon_{n-i_1}\cdots\varepsilon_{n-i_{k_1}}, \qquad Y_2'(n)=
\sum_{\mathbf{i}\in\mathbb{Z}_+^{k_2}}'a^{(2)}(\mathbf{i})
\varepsilon _{n-i_1}\cdots\varepsilon_{n-i_{k_2}},
\end{equation}
where we assume that $a^{(j)}=g^{(j)}L^{(j)}$ as in (\ref{eqa=gL}) is
symmetric, where $g^{(j)}$ is a symmetric GHK(B) with homogeneity
exponent
\[
\alpha_j\in(-k_j/2-1/2,-k_j/2),\qquad j=1,2.
\]
Definition~\ref{DefGHprocess} suggests the following terminology.

\begin{Def}\label{DefHurstExpa^j}
The index
%
\begin{equation}
\label{eqHrange} H=\alpha+k/2+1\in(1/2,1)
\end{equation}
is called the
\emph{associated Hurst index} of the coefficient $a(\cdot)=g(\cdot
)L(\cdot)$ in (\ref{eqa=gL}).
\end{Def}

\begin{Rem}
The associated Hurst indices of the coefficients in $Y_1'(n)$ and
$Y_2'(n)$ will determine the Hurst index of the limit process $Z(t)$ in
(\ref{eqgenerallimit}).
\end{Rem}

We want to obtain a limit theorem for the normalized partial sum of the
product process:
%
\begin{equation}
\label{eqchaosproductproc} X(n):=Y'_1(n)Y_2'(n).
\end{equation}

\begin{Thm}\label{Thmgeneralprodchaos}
Let $X(n)$ be the product process in (\ref{eqchaosproductproc}).
Suppose that $H_j$ is the associated Hurst index of $a^{(j)}(\cdot)$,
$j=1,2$, and assume that $\mathbb{E}|\varepsilon_i|^{4+\delta}<\infty$
for some
$\delta>0$.
\begin{enumerate}[2.]
\item[1.] If $H_1+H_2<3/2$, then $X(n)$ satisfies the
CLT (\ref{eqclt});

\item[2.] If $H_1+H_2>3/2$, then $X(n)$ satisfies
the NCLT (\ref{eqnclt}) with Hurst index $H=H_1+H_2-1$ and limit
%
\begin{equation}
Z(t)=\sum_{r=0}^{k} r!\pmatrix{k_1 \cr r} \pmatrix{k_2\cr r} I_{k_1+k_2-2r}(h_{t,r}),
\end{equation}
where $k=k_1\wedge k_2$ if $k_1\neq k_2$, and $k=k_1-1$ if $k_1=k_2$.
The integrand $h_{t,r}$ above is defined as
%
\begin{equation}
\label{eqhtr}
h_{t,r}(\mathbf{x})= \int_0^t
\bigl( g^{(1)}\otimes_r g^{(2)} \bigr) (s
\mathbf{1}-\mathbf{x})1_{\{s\mathbf{1}>\mathbf{x}\}} \,\mathrm{d}s,
\end{equation}
where
%
\begin{eqnarray}
&& g^{(1)}\otimes_r g^{(2)}(
\mathbf{x})
\nonumber
\\[-8pt]
\label{eqcontractofg}
\\[-8pt]
\nonumber
&&\quad :=\int_{\mathbb{R}_+^r}g^{(1)}(y_1,
\ldots,y_r,x_1,\ldots ,x_{k_1-r})g^{(2)}(y_1,
\ldots,y_r,x_{k_1-r+1},\ldots ,x_{k_2+k_2-2r})\,\mathrm{d}\mathbf{y}
\end{eqnarray}
is a GHK, and when $r=0$, (\ref{eqcontractofg}) is understood as the
tensor product $g^{(1)}\otimes g^{(2)}$. When $r>0$ in (\ref
{eqcontractofg}), we identify $r$ variables of $g^{(1)}$ and
$g^{(2)}$ and integrate over them.
\end{enumerate}
\end{Thm}

This theorem is proved in Section~\ref{secproofchaos}.

\subsection{Limit theorem for a product of long-memory Volterra processes}
Let now
%
\begin{equation}
\label{eqX=product} X(n)=Y_1(n)Y_2(n),
\end{equation}
where
%
\begin{equation}
\label{eqwith-diagonalY} Y_1(n)=\sum_{\mathbf{i}\in\mathbb{Z}_+^{k_1}}a^{(1)}(
\mathbf {i})\varepsilon_{n-i_1}\cdots\varepsilon_{n-i_{k_1}},\qquad
Y_2(n)=\sum_{\mathbf{i}\in\mathbb{Z}_+^{k_2}}a^{(2)}(
\mathbf{i})\varepsilon _{n-i_1}\cdots\varepsilon_{n-i_{k_2}}.
\end{equation}
We assume that $a^{(j)}=g^{(j)}L^{(j)}$ in (\ref{eqa=gL}) is
symmetric, and $g^{(j)}$ is a symmetric GHK(B) with homogeneity
exponent $\alpha_j\in(-k_j/2-1/2,-k_j/2)$, $j=1,2$.
In this case, we can write
\[
X(n)=\sum_{\mathbf{i}\in\mathbb{Z}_+^k}a(\mathbf{i})\varepsilon
_{i_1}\cdots\varepsilon_{i_k},
\]
where $k=k_1+k_2$, and
%
\begin{equation}
\label{eqa=product} a=a^{(1)}\otimes a^{(2)}.
\end{equation}

Let $\mathcal{C}_{1}^2$ to be the collection of partitions of the set
$\{1,\ldots,k_1\}$ such that each set in the partition contains at
least $2$ elements, and similarly let $\mathcal{C}_{2}^2$ be the same
thing for $\{k_1+1,\ldots,k_1+k_2\}$. Any partition $\pi\in\mathcal
{C}_j^2$ can be expressed as $\pi=(P_1,\ldots,P_m)$, where $P_i$,
$i=1,\ldots,m$, are subsets ordered according to their smallest
elements. For example, if $\pi=\{\{1,4\},\{2,3\}\}$, then $P_1=\{1,4\}$
and $P_2=\{2,3\}$.
Let
%
\begin{equation}
\label{eqcsa}
c_j=\sum_{\pi\in\mathcal{C}_j^2} \sum
'_{\mathbf{i}>\mathbf
{0}}a^{(j)}_{\pi}(
\mathbf{i}) \mu_{\pi},\qquad j=1,2,
\end{equation}
where
\[
\mu_\pi=\mu_{p_1}\cdots\mu_{p_m}\qquad\mbox{with }
\mu_p=\mathbb {E}\varepsilon_i^p
\]
and $p_i=|P_i|\ge2$ if $\pi=(P_1,\ldots,P_m)$, and where
$a^{(j)}_{\pi
}(\cdot)$ denotes $a^{(j)}$ with its variables identified according to
the partition $\pi$ (see (\ref{eqapi}) below).

The limit theorem for the normalized partial sum of the centered $X(n)$
in (\ref{eqX=product}) includes several cases. We shall use the
centered multiple Wiener--Stratonovich integral $ \aaa^c_k(\cdot
)$ introduced in (\ref{eqwienerstratonintegral}). The theorem states
that except for some low-dimensional cases (cases 1--4), the limit is up
to some constant the same as the limit for a single factor, namely the
one with the highest $H_j$ (cases 5--7).

\begin{Thm}\label{Thmgeneralprodvolt}
Let $X(n)$ be the product process in (\ref{eqX=product}), where
$a^{(j)}$ has associated Hurst index $H_j=\alpha_j+k_j/2+1\in(1/2,1)$
(Definition~\ref{DefHurstExpa^j}). Assume $\mathbb{E}|\varepsilon
_i|^{2k_1+2k_2+\delta}<\infty$ for some $\delta>0$. Then using the
language of Definition~\ref{Defcltnclt},
\begin{enumerate}[7.]
\item[1.] if $k_1=1$, $k_2=1$, and
$H_1+H_2<3/2$, then $X(n)$ satisfies a CLT (\ref{eqclt});
\item[2.] if $k_1=1$, $k_2=1$, and
$H_1+H_2>3/2$, then $X(n)$ satisfies a NCLT (\ref{eqnclt}) with Hurst
index $H_1+H_2-1$ and limit
\[
Z(t)=\int'_{\mathbb{R}^2}\!\int_0^t
g_1(s-x_1) g_2(s-x_2)
1_{\{
s\mathbf
{1}>\mathbf{x}\}} \,\mathrm{d}s B(\mathrm{d}x_1) B(\mathrm{d}x_2)
\]
(nonsymmetric Rosenblatt process);
\item[3.] if $k_1\ge2$, $k_2=1$, and if $c_1$
in (\ref{eqcsa}) is nonzero, then $X(n)$ satisfies a NCLT (\ref
{eqnclt}) with Hurst index $H_2$ and limit
\[
Z(t)=c_1\int_{\mathbb{R}}\!\int_0^t
g_2(s-x)1_{\{s>x\}} \,\mathrm{d}s B(\mathrm{d}x)
\]
(fractional Brownian motion);

\item[4.] if $k_1= 1$, $k_2\ge2$, and if $c_2$
in (\ref{eqcsa}) is nonzero, then $X(n)$ satisfies a NCLT (\ref
{eqnclt}) with Hurst index $H_1$ and limit
\[
Z(t)=c_2\int_{\mathbb{R}}\!\int_0^t
g_1(s-x)1_{\{s>x\}} \,\mathrm{d}s B(\mathrm{d}x)
\]
(fractional Brownian motion);

\item[5.] if $k_1\ge2$, $k_2\ge2$,
$H_1>H_2$, and if $c_2$ in (\ref{eqcsa}) is nonzero, then $X(n)$
satisfies a NCLT (\ref{eqnclt}) with Hurst index $H_1$, and the\vspace*{-2pt} limit
\[
Z(t)=c_2 \aaa^c_{k_1}(h_{t,1}),
\]
where $h_{t,1}(\mathbf{x})=\int_0^t g_1(s\mathbf{1}-\mathbf{x})
1_{\{
s\mathbf{1}>\mathbf{x}\}} \,\mathrm{d}s$;
\item[6.] if $k_1\ge2$, $k_2\ge2$,
$H_1<H_2$, and if $c_1$ in (\ref{eqcsa}) is nonzero, then $X(n)$
satisfies a NCLT (\ref{eqnclt}) with Hurst index $H_2$, and the limit
\[
Z(t)=c_1 \aaa^c_{k_2}(h_{t,2}),
\]
where $h_{t,2}(\mathbf{x})=\int_0^t g_2(s\mathbf{1}-\mathbf{x})
1_{\{
s\mathbf{1}>\mathbf{x}\}} \,\mathrm{d}s$;

\item[7.] if $k_1\ge2$, $k_2\ge2$,
$H_1=H_2$, and if at least one of the $c_j$'s in (\ref{eqcsa}) is
nonzero, then $X(n)$ satisfies a NCLT (\ref{eqnclt}) with Hurst index
$H_1=H_2$, and the limit
\[
Z(t)=c_1 \aaa^c_{k_2}(h_{t,2})+c_2
\aaa^c_{k_1}(h_{t,1}).
\]
\end{enumerate}
\end{Thm}

\begin{Rem}
These constants $c_j$'s in the theorem are nonzero if, for example,
every $a^{(j)}(\mathbf{i})>0$, $j=1,2$.
\end{Rem}

The theorem, which is proved in Section~\ref{secproofvolt}, seems
bewildering at first glance. But there is structure into it. The cases~3 and 4 are symmetric,
and so are the cases~5 and~6. Case~1 involves
short-range dependence, while all the other cases involve long-range
dependence. Case~2 involves the nonsymmetric
Rosenblatt process, originally introduced by
Maejima and Tudor
\cite{maejimatudor2012selfsimilar}. Cases~3 and
4 involve fractional Brownian motion since one
of the orders $k$ equals $1$. The typical cases are 5 (and~6). In these cases,
quite surprisingly, it is not the orders $k_1$ or $k_2$ that matter,
but the process $Y_1(n)$ or $Y_2(n)$ in (\ref{eqX=product}) with the
highest value of $H$. In the boundary case~7, where $H_1=H_2$, they both contribute.\vspace*{-3pt}

\subsection{Limit theorem for the mixed case}
Now we consider the mixed case (\ref{eqprodofmixed}), where
$Y_1'(n)$ is as in (\ref{eqoff-diagonalY}) and $Y_2(n)$ is as in
(\ref
{eqwith-diagonalY}).
Let
%
\begin{equation}
\label{eqX=productmixed} X(n)=Y_1'(n)Y_2(n).
\end{equation}

We only state the case which does not overlap Theorem~\ref{Thmgeneralprodchaos} and Theorem~\ref{Thmgeneralprodvolt}, that is, both
$Y_1'(n)$ and $Y_2(n)$ are nonlinear: $k_1\ge2$ and $k_2\ge2$. The
limit, up to some constant, turns out to be the same as the limit for
the single factor $Y_1'(n)$.

\begin{Thm}\label{Thmgeneralprodmixed}
Let $X(n)$ be the product process in (\ref{eqX=productmixed}), where
$a^{(j)}$ has associated Hurst index $H_j=\alpha_j+k_j/2+1\in(1/2,1)$
(Definition~\ref{DefHurstExpa^j}). Assume $k_1\ge2$, $k_2\ge2$
and $\mathbb{E}|\varepsilon_i|^{2+2k_2+\delta}<\infty$ for some
$\delta>0$. Then
using the language of Definition~\ref{Defcltnclt}, if $c_2$ in (\ref
{eqcsa}) is nonzero, then $X(n)$ satisfies a NCLT (\ref{eqnclt})
with Hurst index $H_1=\alpha_1+k_1/2+1$, and the limit is
\[
Z(t)=c_2 I_{k_1}(h_{t,1}),
\]
where\vadjust{\goodbreak} $h_{t,1}(\mathbf{x})=\int_0^t g_1(s\mathbf{1}-\mathbf{x})
1_{\{
s\mathbf{1}>\mathbf{x}\}} \,\mathrm{d}s$.
\end{Thm}

Theorem~\ref{Thmgeneralprodmixed} is proved in Section~\ref{secproofmixed}.

\begin{Rem}
If the noises $\varepsilon_i$'s are Gaussian, then the normalized partial
sum
\[
\frac{1}{A(N)}\sum_{n=1}^{[Nt]}X(n)
\]
considered in Theorem~\ref{Thmgeneralprodchaos} \ref{Thmgeneralprodvolt} and \ref{Thmgeneralprodmixed} belongs to a Wiener chaos
of finite order. There is a rich literature on obtaining Berry--Esseen
type quantitative limit theorems for elements on Wiener chaos. For the
case where the limit is Gaussian, see the monograph \cite{nourdinpeccati2012normal} and the references therein; for the case
where the limit belongs to higher-order Wiener chaos, see
\cite{davydovmartynova1987limit,breton2006convergence} and
\cite{nourdinpoly2013convergence}. The case where $\varepsilon_i$'s
are non-Gaussian may also be treated using techniques from \cite{nourdinpeccati2010invariance}.

The quantitative results mentioned above, however, seem not directly
applicable to the limit theorems considered here. This is because, as
it will be clear in the proofs of these theorems, $\frac{1}{A(N)}\sum_{n=1}^{[Nt]}X(n)$ does not have a ``clean'' structure as that
considered in the works mentioned above. In particular, the
decomposition of $\frac{1}{A(N)}\sum_{n=1}^{[Nt]}X(n)$ yields many
terms. Some of the quantitative results mentioned above may be
applicable to the terms which contribute to the limit, but there are
other terms in the decomposition which converge in $L^2(\Omega)$ to
zero. How to deal with these degenerate terms is an open problem.
\end{Rem}

\section{Preliminary results}\label{secprelim}
A central idea in establishing the limit theorems is to involve the
\emph{nonsymmetric discrete chaos process} which generalizes the chaos
process in (\ref{eqXnintro1}) by allowing different sequences of
noises. We shall now define it.
Let $\bolds{\varepsilon}_i=(\varepsilon_i^{(1)},\ldots,\varepsilon
_i^{(k)})$ be
an i.i.d. vector where each component has mean $0$ and finite
variance. The components $\varepsilon_i^{(1)},\ldots,\varepsilon_i^{(k)}$ are
typically dependent. Introduce the following nonsymmetric discrete
chaos process
%
\begin{equation}
\label{eqXnnonsym} Y'(n)=\sum_{0<i_1,\ldots,i_k<\infty}'
a(i_1,\ldots,i_k)\varepsilon ^{(1)}_{n-i_1}
\cdots\varepsilon^{(k)}_{n-i_k},
\end{equation}
where $\sum_{\mathbf{i}\in\mathbb{Z}_+^k}' a(\mathbf{i})^2<\infty
$ so
that $X'(n)$ is well-defined in the $L^2(\Omega)$-sense. Let
\[
\Sigma(i,j)= \mathbb{E}\varepsilon_n^{(i)}
\varepsilon_n^{(j)}.
\]
The autocovariance of $Y'(n)$ is then given by
%
\begin{equation}
\label{eqXnnonsymautocov} \hspace*{-6pt}\gamma(n)=\sum_{\sigma}\sum
_{0<i_1,\ldots,i_k<\infty}' a(i_1,\ldots
,i_k) a(i_{\sigma(1)}+n,\ldots,i_{\sigma(k)}+n) \Sigma
(i_1,i_{\sigma
(1)})\cdots\Sigma(i_k,i_{\sigma(k)}),
\end{equation}
where in the summation $\sigma$ runs over all the $k!$ permutations of
$\{1,\ldots,k\}$.
The following lemma is useful for studying the asymptotic properties of
the covariance of $X'(n)$.

\begin{Lem}\label{Lembound}
Suppose that in (\ref{eqXnnonsym}), there exist constant $c_0>0$
and $\gamma_j<-1/2$, $j=1,\ldots,k$, such that
%
\begin{equation}
\label{eqabound} \bigl\llvert a(i_1,\ldots,i_k)\bigr\rrvert
\le c_0 i_1^{\gamma_1}\cdots i_k^{\gamma_k}.
\end{equation}
Let
%
\begin{equation}
\label{eqH^*} H^*=\alpha+k/2+1 \qquad\mbox{with }\alpha=\sum
_{j=1}^{k}\gamma_j.
\end{equation}
\begin{itemize}
\item If $H^*<1/2$, then $\sum_{n=-\infty}^\infty|\gamma(n)|<\infty$,
and $\operatorname{Var}[\sum_{n=1}^N Y'(n)]\le c_1 N$ for some $c_1>0$;
\item If $H^*>1/2$, then $|\gamma(n)|\le c_2 n^{2H^*-2}$ for some
$c_2>0$, and $\operatorname{Var}[\sum_{n=1}^N Y'(n)]\le c_3 N^{2H^*}$
for some $c_3>0$.
\end{itemize}
\end{Lem}

\begin{pf}
The case $H^*<1/2$ was proved in Proposition~5.4 in \cite{baitaqqu2014convergence}.

In the case $H^*>1/2$, let $\widetilde{|a|}$ be the symmetrization of
$|a|(\mathbf{i}):=|a(\mathbf{i})|$, then for $n\ge0$, by (\ref
{eqXnnonsymautocov}) and (\ref{eqabound}),
\begin{eqnarray*}
\bigl\llvert \gamma(n)\bigr\rrvert &\le& C_0 \sum
_{\mathbf{i}\in\mathbb{Z}_+^k} \widetilde{|a|}(\mathbf{i}+n\mathbf{1})
\widetilde{|a|}(\mathbf{i})
\\
&\le& C_1 \sum_{\sigma} \sum
_{i_1=1}^\infty\cdots\sum_{i_k=1}^\infty
(i_1+n)^{\gamma_1}\cdots(i_k+n)^{\gamma_k}
i_1^{\gamma_{\sigma(1)}} \cdots i_k^{\gamma_{\sigma(k)}}
\\
&\le& C_2 \sum_{\sigma} n^{\gamma_1+\gamma_{\sigma(1)}+1}
\cdots n^{\gamma
_k+\gamma_{\sigma(k)}+1}=C_3 n^{2\alpha+k}=C_3n^{2H^*-2},
\end{eqnarray*}
where\vspace*{1pt} the $C_i$'s are positive constants, and $\sigma$ in the summation
runs over all the permutations of $\{1,\ldots,k\}$. $\operatorname
{Var}[\sum_{n=1}^N
Y'(n)]\le c_3 N^{2H^*}$ then follows as a standard result.
\end{pf}

\begin{Rem}
In the applications of Lemma~\ref{Lembound}, the inequality (\ref{eqabound}) is often not seen in this form. For example, the function
$a(\cdot)$ defined on $\mathbb{Z}_+^k$ may satisfy
\[
\bigl\llvert a(\mathbf{i})\bigr\rrvert \le C (i_1+
\cdots+i_{k_1})^{\alpha
_1}(i_{k_1+1}+\cdots
+i_{k_1+k_2})^{\alpha_2},
\]
for some $C>0$,
where $k_1+k_2=k$, and $\frac{\alpha_j}{k_j}<-\frac{1}{2}$, then it is
easily verified by the arithmetic--geometric mean inequality
\[
k^{-1}\sum_{j=1}^k
y_j \ge \Biggl(\prod_{j=1}^k
y_j \Biggr)^{1/k}
\]
for $y_j>0$, that (\ref{eqabound}) is satisfied since $\alpha<0$. It
is also verified for a function $a_\pi(\cdot)$ which is $a(\cdot)$ with
some of its variables identified.

In general when applying Lemma~\ref{Lembound}, we will omit the
verification of (\ref{eqabound}) which usually can be easily done as
indicated above. We will merely count the \emph{total homogeneity
exponents} of the bound, which in the preceding example is $\alpha
=\alpha_1+\alpha_2$.
\end{Rem}

For convenience, we make the following definition.

\begin{Def}\label{DefHurstexpXn}
Let $X(n)$ be a stationary process with mean $0$ and finite variance.
We say
\begin{itemize}
\item
$X(n)$ has a memory parameter of \emph{at most} (denoted using $\le$)
$H$, if
\[
\operatorname{Var} \Biggl[\sum_{n=1}^N
X(n) \Biggr]\le c N^{2H}
\]
for some $c>0$;
\item
$X(n)$ has a memory parameter (denoted using $=$) $H$, if
\[
\operatorname{Var} \Biggl[\sum_{n=1}^N
X(n) \Biggr]\sim c N^{2H}
\]
as $N\rightarrow\infty$ for some $c>0$.
\end{itemize}
\end{Def}

\begin{Rem}
In view of the definition above, Lemma~\ref{Lembound} states that if
$Y'(n)$ in (\ref{eqXnintro1}) satisfies (\ref{eqabound}), then
$Y'(n)$ has a memory parameter of at most $1/2$ if $H^*<1/2$ and of at
most $H^*$ if $H^*>1/2$.
\end{Rem}

\begin{Pro}[(Proposition~5.4 of \cite{baitaqqu2014convergence})]\label
{ThmcltXn}
Let $Y'(n)$ be given as in (\ref{eqXnnonsym}) with coefficient
$a(\cdot)$ satisfying (\ref{eqabound}) and $H^*< 1/2$ in Lemma~\ref
{Lembound}. Then
\[
N^{-1/2}\sum_{n=1}^{[Nt]}
\bigl[Y'(n)-\mathbb{E}Y'(n)\bigr]\stackrel {f.d.d.} {
\longrightarrow}\sigma B(t),
\]
where
\[
\sigma^2=\sum_{n=-\infty}^\infty
\operatorname{Cov}\bigl[Y'(n),Y'(0)\bigr],
\]
$B(t)$ is a standard Brownian motion, and $\stackrel
{f.d.d.}{\longrightarrow}$ stands for
convergence of finite-dimensional distributions.

If each $\varepsilon_i^{(1)},\ldots,\varepsilon_{i}^{(k)}$ has a moment
greater than $2$, then the tightness of
\[
N^{-1/2}\sum_{n=1}^{[Nt]}
\bigl[Y'(n)-\mathbb{E}Y'(n)\bigr]
\]
in $D[0,1]$ holds and thus $\stackrel{f.d.d.}{\longrightarrow}$ can
be replaced by weak
convergence $\Rightarrow$ in $D[0,1]$.
\end{Pro}

\begin{Rem}
The above $\stackrel{f.d.d.}{\longrightarrow}$ or $\Rightarrow$
convergence also holds for a
linear combination of different $Y'(n)$'s defined on a common i.i.d.
noise vector $\bolds{\varepsilon}_i$, while the $Y'(n)$'s can have
different orders and involve different subvectors of $\bolds{\varepsilon
}_i$, provided the coefficient of each $Y'(n)$ satisfies (\ref{eqabound}) with $H^*<1/2$.
\end{Rem}

We now state an important result concerning the weak convergence of a
discrete chaos to a Wiener chaos. Let $h$ be a function defined on
$\mathbb{Z}^k$ such that $\sum'_{\mathbf{i}\in\mathbb{Z}^k_+}
h(\mathbf
{i})^2<\infty$, where $'$ indicates the exclusion of the diagonals
$i_p=i_q$, $p\neq q$.
Let $Q_k(h)$ be defined as follows:
%
\begin{eqnarray}
\label{eqQkh} Q_k(h)=Q_k(h,\bolds{\varepsilon})=\sum
'_{(i_1,\ldots,i_k)\in\mathbb{Z}^k} h(i_1,
\ldots,i_k) \varepsilon_{i_1}\cdots\varepsilon_{i_k}=
\sum'_{\mathbf
{i}\in\mathbb{Z}^k} h(\mathbf{i})\prod
_{p=1}^k \varepsilon_{i_p},
\end{eqnarray}
where $\varepsilon_i$'s are i.i.d. noises. Observe that $Q_k(h)$ is
invariant under permutation of the arguments of $h(i_1,\ldots,i_k)$. So
if $\tilde{h}$ is the symmetrization of $h$, then $Q_k(h)=Q_k(\tilde{h})$.

Suppose now that we have a sequence of function vectors $\mathbf
{h}_n=(h_{1,n},\ldots,h_{J,n})$ where each $h_{j,n}\in L^2(\mathbb
{Z}^{k_j})$, $j=1,\ldots,J$.

\begin{Pro}[(Proposition~4.1 of \cite{baitaqqu2013generalized})]\label
{ProPoly->Wiener}
Let
\[
\tilde{h}_{j,n}(\mathbf{x})=n^{k_j/2}h_{j,n} \bigl([n
\mathbf {x}]+\mathbf {c}_j \bigr),\qquad j=1,\ldots,J,
\]
where $\mathbf{c}_j\in\mathbb{Z}^k$. Suppose that there exists
$h_j\in
L^2(\mathbb{R}^{k_j})$, such that
%
\begin{equation}
\label{eqhtildehL2conv} \llVert \tilde{h}_{j,n}-h_j\rrVert
_{L^2(\mathbb{R}^{k_j})}\rightarrow0
\end{equation}
as $n\rightarrow\infty$. Then, as $n\rightarrow\infty$, we have the
following joint convergence in distribution:
\begin{eqnarray*}
&&\mathbf{Q}:= \bigl(Q_{k_1}(h_{1,n}),\ldots,Q_{k_J}(h_{J,n})
\bigr) \stackrel{d} {\rightarrow} \mathbf{I}:= \bigl(I_{k_1}(h_1),
\ldots,I_{k_J}(h_J) \bigr).
\end{eqnarray*}
\end{Pro}

\section{Proofs}\label{secproof}
\subsection{Proof of Theorem \texorpdfstring{\protect\ref{Thmgeneralprodchaos}}{3.5} where
diagonals are excluded}\label{secproofchaos}

We first show that $g^{(1)}\otimes_r g^{(2)}$ in (\ref{eqcontractofg}) is a GHK.

\begin{Lem}\label{LemgcontractgisGHK}
Let $g^{(j)}$ be a symmetric GHK(B) with homogeneity exponent $\alpha
_j$ defined on $\mathbb{R}_+^{k_j}$, $j=1,2$. Suppose in addition that
either $k_1\ge2$ or $k_2\ge2$, and that
%
\begin{equation}
\label{eqalpha1+alpha2>} \alpha_1+\alpha_2>-(k_1+k_2+1)/2,
\end{equation}
and set
\[
r= %
\cases{ 0,\ldots,k_1\wedge k_2, &\quad $
\mbox{if } k_1\neq k_2$,\vspace*{3pt}
\cr
0,
\ldots,k_1-1, & \quad$\mbox{if } k_1= k_2$.}
\]
If the function $g^{(1)}\otimes_r g^{(2)} $
is nonzero, then it is a GHK on $\mathbb{R}_+^{k_1+k_2-2r}$ with
homogeneity exponent $\alpha_1+\alpha_2+r$.
\end{Lem}

\begin{pf}
When $r=0$, $g^{(1)}\otimes g^{(2)}$ is a tensor product of two
GHK(B)s. It is a GHK because condition~1 of Definition~\ref{DefGHK} is
satisfied with homogeneity\vspace*{-2pt} exponent
%
\begin{equation}
\label{eqalpha1+alpha2+k+1} -(k_1+k_2+1)/2<\alpha_1+
\alpha_2<-(k_1+k_2)/2
\end{equation}
[see (\ref{eqalpha})], and condition~2 of Definition~\ref{DefGHK} is
satisfied because
\begin{eqnarray*}
&&\int_{\mathbb{R}_+^{k_1+k_2}}\bigl\llvert g^{(1)}(
\mathbf{x}_1)g^{(2)}(\mathbf {x}_2)g^{(1)}(
\mathbf{1}+\mathbf{x}_1) g^{(2)}(\mathbf{1}+\mathbf
{x}_2)\bigr\rrvert \,\mathrm{d}\mathbf{x}_1\,\mathrm{d}\mathbf{x}_2
\\[-2pt]
&&\quad=\int_{\mathbb{R}_+^{k_1}} \bigl\llvert g^{(1)}(
\mathbf{x})g^{(1)}(\mathbf {1}+\mathbf{x})\bigr\rrvert \,\mathrm{d}\mathbf{x} \int
_{\mathbb{R}_+^{k_2}} \bigl\llvert g^{(2)}(\mathbf
{x})g^{(2)}(\mathbf{1}+\mathbf{x})\bigr\rrvert \,\mathrm{d}\mathbf{x}<\infty.
\end{eqnarray*}
We shall now focus on the case $r>0$.

Consider\vspace*{1pt} first $k_1\ge2$ and $k_2=1$ (the case $k_1=1$ and $k_2\ge2$
is similar), so that $g^{(2)}(x)=C x^{\alpha_2}$ for some $C\neq0$,
where $\alpha_2\in(-1,-1/2)$. Fix an $\mathbf{x}=(x_1,\ldots
,x_{k-1})\in\mathbb{R}_+^{k_1-1}$, then
\[
\int_0^\infty\bigl\llvert g^{(1)}(y,
\mathbf{x})\bigr\rrvert y^{\alpha_2} \,\mathrm{d}y\le C \int_0^\infty
(y+x_1\cdots+x_{k_1-1})^{\alpha_1}y^{\alpha_2} \,\mathrm{d}y <
\infty,
\]
because near $y=0$ (the other $\mathbf{x}>\mathbf{0}$), the integrand
behaves like $y^{\alpha_2}$, where $\alpha_2>-1$, while near
$y=\infty
$, the integrand is like $y^{\alpha_1+\alpha_2}$, where $\alpha_1<-1$
and $\alpha_2<-1/2$. Hence, $g^{(1)}\otimes_1 g^{(2)}$ is well-defined
in this case. It is easy to check that
\[
g^{(1)}\otimes_1 g^{(2)}(\lambda\mathbf{x})=
\lambda^{\alpha
_1+\alpha
_2+1} g^{(1)}\otimes_1 g^{(2)}(
\mathbf{x})
\]
for any $\lambda>0$ by using a change of variable and using the
homogeneity of $g^{(j)}$. We are left to show that $g:=g^{(1)}\otimes_1
g^{(2)}$ satisfies condition~2 of Definition~\ref{DefGHK}.
This is true because the function $f(x):=\int_0^\infty(x+y)^{\alpha_1}
y^{\alpha_2}\,\mathrm{d}y$ is $f(x)=C_0 x^{\alpha_1+\alpha_2+1}$ for some
$C_0>0$. So
%
\begin{equation}
\label{eqboundoneleft} \bigl\llvert g^{(1)}\otimes_1
g^{(2)}(\mathbf{x})\bigr\rrvert \le C (x_1+\cdots
+x_{k_1-1})^{\alpha_1+\alpha_2+1}=:g^*(\mathbf{x})
\end{equation}
for some $C>0$. Note that $g^*(\cdot)$ is a GHK(B) on $\mathbb
{R}^{k_1-1}$ with
\[
-(k_1-1)/2-1/2<\alpha_1+\alpha_2+1<-(k_1-1)/2
\]
because $\alpha_1<-1/2$, $\alpha_2<-k_2/2$ and $\alpha_1+\alpha
_2>-(1+k_2+1)/2$ by assumption (\ref{eqalpha1+alpha2>}). So
$g=g^{(1)}\otimes_1 g^{(2)}$ satisfies condition 2 of Definition~\ref
{DefGHK} because the dominating function $g^*$ does.

Suppose now that $k_1\ge2$ and $k_2\ge2$. Consider first the case $1
\le r\le(k_1\wedge k_2)-1 $. Using the bound $g^{(j)}(\mathbf{x})\le
C\|\mathbf{x}\|^{\alpha_j}$, one has by applying Cauchy--Schwarz and
integrating power functions iteratively that
\begin{eqnarray}
&&\bigl\llvert g^{(1)}\otimes_r g^{(2)}(\mathbf{x})
\bigr\rrvert
\nonumber
\\
&&\quad\le C \int_{\mathbb{R}_+^r}(y_1+\cdots+y_r+x_1+
\cdots +x_{k_1-r})^{\alpha_1}
\nonumber
\\
\label{eqitercauchyschwartz}
&&\hspace*{24pt}\qquad{}\times(y_1+\cdots +y_r+x_{k_1-r+1}+\cdots
+x_{k_2+k_2-2r})^{\alpha_2}\,\mathrm{d}y_1\cdots \,\mathrm{d}y_r
\\
\nonumber
&&\quad= C \int_{\mathbb{R}^{r-1}_+}\,\mathrm{d}y_1\cdots
\,\mathrm{d}y_{r-1} \biggl(\int_0^\infty
(y_1+\cdots+y_r+x_1+\cdots+x_{k_1-r})^{2\alpha_1}\,\mathrm{d}y_r
\biggr)^{1/2}
\nonumber
\\
&&\hspace*{32pt}\qquad{}\times \biggl(\int_0^\infty(y_1+
\cdots+y_r+x_{k_1-r+1}+\cdots +x_{k_2+k_2-2r})^{2\alpha_2}\,\mathrm{d}y_r
\biggr)^{1/2}
\nonumber
\\
&&\quad\le C \int_{\mathbb{R}_+^{r-1}}(y_1+\cdots+y_{r-1}+x_1+
\cdots +x_{k_1-r})^{\alpha_1+1/2}
\nonumber
\\
&&\hspace*{30pt}\qquad{}\times(y_1+\cdots+y_{r-1}+x_{k_1-r+1}+\cdots
+x_{k_2+k_2-2r})^{\alpha
_2+1/2}\,\mathrm{d}\mathbf{y}
\nonumber
\\
&&\qquad \cdots
\nonumber
\\
\nonumber
&&\quad\le C (x_1+\cdots+x_{k_1-r})^{\alpha_1+r/2}(x_{k_1-r+1}+
\cdots +x_{k_1+k_2-2r})^{\alpha_2+r/2}=:g^*(\mathbf{x}).
\end{eqnarray}
The dominating function $g^*$ is a GHK because it is a tensor product
of two GHK(B)'s on $\mathbb{R}_+^{k_j}$, $j=1,2$, and
\[
-\frac{(k_1-r)+(k_2-r)+1}{2}<(\alpha_1+r/2)+(\alpha_2+r/2)<-
\frac
{(k_1-r)+(k_2-r)}{2},
\]
as in the inequality (\ref{eqalpha1+alpha2+k+1}).
Therefore, the bound $g^*(\mathbf{x})$, and hence
the kernel $g^{(1)}\otimes_r g^{(2)}$ satisfy condition 2 of Definition~\ref{DefGHK}. Moreover, the homogeneity exponent of $g^{(1)}\otimes_r
g^{(2)}$ is $\alpha_1+\alpha_2+r$ in condition~1 of Definition~\ref
{DefGHK}. This can be easily verified as above by change of variables
and using the homogeneity of $g^{(j)}$.

The only case left is: $k_1\neq k_2 \ge2$ and $r=k_1\wedge k_2$.
Suppose $k_1<k_2$. In this case, condition~2 of Definition~\ref
{DefGHK} can be checked by first applying the iterative
Cauchy--Schwarz argument leading to (\ref{eqitercauchyschwartz})
until only one variable of $g^{(1)}$ is unintegrated, and then bounding
the last fold of integration similarly as in (\ref{eqboundoneleft}).
Hence, in this case as well, $g^{(1)}\otimes_r g^{(2)}$ is GHK.
\end{pf}

The following lemma shows a noncentral convergence involving
$g^{(1)}\otimes_r g^{(2)}$ appearing in (\ref{eqcontractofg}).

\begin{Lem}\label{Lemchaosproductnclt}
Suppose that all the assumptions in Lemma~\ref{LemgcontractgisGHK}
hold. Let $a^j(\cdot)=g^{(j)}L^{(j)}$, $j=1,2$, be as assumed before. Set
\begin{eqnarray*}
X'_r(n)&:= &\sum_{(\mathbf{u},\mathbf{i})>\mathbf{0}}'
a^{(1)}(u_1,\ldots ,u_r,i_1,
\ldots,i_{k_1-r})
\\
&&\hspace*{4pt}\qquad{}\times a^{(2)}(u_1,\ldots ,u_r,i_{k_1-r+1},
\ldots,i_{k_1+k_2-2r})\varepsilon_{n-i_1}\cdots \varepsilon _{n-i_{k_1+k_2-2r}},
\end{eqnarray*}
where $\varepsilon_i$'s are i.i.d. with mean $0$ and variance $1$. We
then have
\[
\frac{1}{N^{H}}\sum_{n=1}^{[Nt]}
X_r'(n) \stackrel {f.d.d.} {\longrightarrow}Z_r(t):=
I_{k_1+k_2-2r}(h_{t,r})
\]
jointly for all the $r=0,1\ldots,k$ where $k$ is as defined in Theorem~\ref{Thmgeneralprodchaos}, and
where
\[
H=\alpha_1+\alpha_2+(k_1+k_2)/2+1
\in(1/2,1).
\]
\end{Lem}

\begin{pf}
In view of Proposition~\ref{ProPoly->Wiener}, we need only to prove
the convergence for a single $r$ and a single $t>0$, and the joint
convergence for different $r$'s and $t$'s follows.
We assume for simplicity that $a^{(j)}(\cdot)=g^{(j)}(\cdot)$ (setting
$L=1$), and including a general $L$ in (\ref{eqa=gL}) is easy.
We focus on the case $r\ge1$, since the case $r=0$ follows from
Theorem~6.5 of \cite{baitaqqu2013generalized}, although the proof
for case $r=0$ may be regarded as contained in the proof below with
$\mathbf{u}$ being an empty vector.

Let $\mathbf{u}=(u_1,\ldots,u_r)$, $\mathbf{i}_1=(i_1,\ldots
,i_{k_1-r})$, $\mathbf{i}_2=(i_{k_1-r+1},\ldots,i_{k_1+k_2-2r})$, and
$\mathbf{i}=(\mathbf{i}_1,\mathbf{i}_2)$.
We define the sum
\[
\sum_{n=1}^{Nt}x_n:=\sum
_{n=1}^{[Nt]}x_n + \bigl(Nt-[Nt]
\bigr)x_{[Nt]+1}=N\int_0^t
x_{1+[Ny]}\,\mathrm{d}y.
\]
Obviously,
\[
\mathbb{E} \Biggl[\frac{1}{N^{H}}\sum_{n=1}^{[Nt]}
X_r'(n)-\frac
{1}{N^{H}}\sum
_{n=1}^{Nt} X_r'(n)
\Biggr]^2\rightarrow0
\]
as $N\rightarrow\infty$. One can thus focus on $\frac{1}{N^{H}}\sum_{n=1}^{Nt} X_r'(n)$ instead.
\begin{eqnarray*}
\frac{1}{N^{H}}\sum_{n=1}^{Nt}
X_r'(n)&=&\sum_{\mathbf{i}\in\mathbb
{Z}^{k_1+k_2-2r} }'
\frac{1}{N^{H}}\sum_{n=1}^{Nt} \sum
_{\mathbf
{u}\in
D(\mathbf{i},n)} g^{(1)}(\mathbf{u},n\mathbf{1}-
\mathbf{i}_1)1_{\{
n\mathbf{1}>\mathbf{i}_1\}}
\\
&&\hspace*{86pt}\qquad{}\times g^{(2)}(\mathbf{u},n\mathbf {1}-\mathbf{i}_2)1_{\{n\mathbf{1}>\mathbf{i}_2\}}
\prod_{j=1}^{k_1+k_2-r}\varepsilon_{i_j}
\\
&=:&Q_{k_1+k_2-2r}(h_{N,t,r}),
\end{eqnarray*}
using the notation (\ref{eqQkh}), where
\[
h_{N,t,r}(\mathbf{i}):=\frac{1}{N^{H}}\sum
_{n=1}^{Nt} \sum_{\mathbf
{u}\in D(\mathbf{i},n)}
g^{(1)}(\mathbf{u},n\mathbf{1}-\mathbf {i}_1)g^{(2)}(
\mathbf{u},n\mathbf{1}-\mathbf{i}_2)1_{\{n\mathbf
{1}>\mathbf{i}\}}
\]
and
\[
D(\mathbf{i},n)=\bigl\{\mathbf{u}\in\mathbb{Z}_+^r\dvt u_p
\neq u_q \mbox{ if }p\neq q; \mbox{ and } u_p\neq n-
i_q \mbox{ even if } p=q\bigr\}.
\]
Set $\mathbf{x}_1\in\mathbb{R}^{k_1-r}$, $\mathbf{x}_1\in\mathbb
{R}^{k_2-r}$ and $\mathbf{x}=(\mathbf{x}_1,\mathbf{x}_2)$. Define
\[
E(\mathbf{x},N)=\bigl\{\mathbf{u}\in\mathbb{Z}_+^r\dvt u_p
\neq u_q \mbox{ if }p\neq q; \mbox{ and } u_p\neq n-
[Nx_q]-1 \mbox{ even if } p=q\bigr\}.
\]
In view of Proposition~\ref{ProPoly->Wiener} and using the homogeneity
of $g^{(j)}$'s, one writes:
\begin{eqnarray*}
\tilde{h}_{N,t,r}(\mathbf{x})&=&N^{(k_1+k_2-2r)/2}h_{N,t}
\bigl([N\mathbf {x}]+\mathbf{1} \bigr)
\\
&=&\frac{1}{N^{\alpha_1+\alpha_2+ r+1}}\sum_{n=1}^{Nt} \sum
_{\mathbf
{u}\in E(\mathbf{x},n)} g^{(1)}\bigl(\mathbf{u},n
\mathbf{1}-[N\mathbf {x}_1]-\mathbf{1}\bigr)
g^{(2)}\bigl(\mathbf{u},n\mathbf{1}-[N\mathbf
{x}_2]-\mathbf{1}\bigr)1_{\{n\mathbf{1}>\mathbf{i}\}}
\\
&=&\sum_{n=1}^{Nt} \frac{1}{N} \sum
_{\mathbf{u}\in E(\mathbf{x},n)} \frac
{1}{N^r} g^{(1)} \biggl(
\frac{\mathbf{u}}{N},\frac{n\mathbf
{1}-[N\mathbf
{x}_1]-\mathbf{1}}{N} \biggr)
g^{(2)} \biggl(\frac{\mathbf
{u}}{N},\frac{n\mathbf{1}-[N\mathbf{x}_2]-\mathbf{1}}{N}
\biggr)1_{\{
n\mathbf{1}>\mathbf{i}\}}
\\
&=&\int_0^t \mathrm{d}s \int_{\mathbb{R}_+^r}\,\mathrm{d}
\mathbf{y} g^{(1)} \biggl(\frac
{[N\mathbf{y}]+\mathbf{1}}{N},\frac{[Ns]\mathbf{1}-[N\mathbf
{x}_1]}{N} \biggr)
\\
&&\hspace*{43pt}{}\times g^{(2)} \biggl(\frac{[N\mathbf
{y}]+\mathbf
{1}}{N},\frac{[Ns]\mathbf{1}-[N\mathbf{x}_2]}{N}
\biggr)1_{\{
[Ns]\mathbf
{1}>[N\mathbf{x}]\}\cap F(N)},
\end{eqnarray*}
where we correspond $\mathbf{u}$ to $[N\mathbf{y}]+\mathbf{1}$, $n$ to
$[Ns]+1$, and
\begin{eqnarray*}
F(N)&=&\bigl\{( \mathbf{x},\mathbf{y},s)\dvt [Ny_p]
\neq[Ny_q], [Nx_p]\neq [Nx_q],
\\
&&\hspace*{4pt} \mbox{ if }p\neq q; \mbox{ and } [Ny_p]\neq[Ns]-[Nx_q] \mbox{ even if }p=q
\bigr\}.
\end{eqnarray*}
In view of Proposition~\ref{ProPoly->Wiener}, the goal is to show that
%
\begin{equation}
\label{eqgoalcontractL^2} \lim_{N\rightarrow\infty}\llVert \tilde{h}_{N,t,r}-h_{t,r}
\rrVert _{L^2(\mathbb
{R}^{k_1+k_2-2r})}=0,
\end{equation}
where $h_{t,r}$ is given in (\ref{eqhtr}).
By the a.e. continuity of $g^{(j)}$'s and the fact that
$1_{F(N)}\rightarrow1$ a.e. as $N\rightarrow\infty$, one has
\begin{eqnarray*}
&&g^{(1)} \biggl(\frac{[N\mathbf{y}]+\mathbf{1}}{N},\frac{[Ns]\mathbf
{1}-[N\mathbf{x}_1]}{N}
\biggr)g^{(2)} \biggl(\frac{[N\mathbf
{y}]+\mathbf
{1}}{N},\frac{[Ns]\mathbf{1}-[N\mathbf{x}_2]}{N}
\biggr)1_{\{
[Ns]\mathbf
{1}>[N\mathbf{x}]\}\cap F(N)}
\\
&&\quad \rightarrow g^{(1)} (\mathbf{y},s\mathbf{1}-\mathbf{x}_1
)g^{(2)} (\mathbf{y},s\mathbf{1}-\mathbf{x}_2
)1_{\{s\mathbf
{1}>\mathbf{x}\}
}\qquad \mbox{for a.e. } (\mathbf{x},\mathbf{y},s).
\end{eqnarray*}
We are left to establish suitable bound to apply the dominated
convergence theorem.
To this end, since $g^{(j)}(\mathbf{x})\le C \|\mathbf{x}\|^{\alpha
_j}=:g^{(j)*}(\mathbf{x})$ on $\mathbb{R}_+^{k_j}$, we have the
following bound:
%
\begin{eqnarray}
&&\biggl\llvert g^{(1)} \biggl(\frac{[N\mathbf{y}]+\mathbf{1}}{N},\frac
{[Ns]\mathbf
{1}-[N\mathbf{x}_1]}{N}
\biggr)g^{(2)} \biggl(\frac{[N\mathbf
{y}]+\mathbf
{1}}{N},\frac{[Ns]\mathbf{1}-[N\mathbf{x}_2]}{N} \biggr)\biggr
\rrvert 1_{\{
[Ns]\mathbf{1}>[N\mathbf{x}]\}\cap F(N)}
\nonumber
\\
&&\quad\le g^{(1)*} \biggl(\frac{[N\mathbf{y}]+\mathbf{1}}{N},\frac
{[Ns]\mathbf
{1}-[N\mathbf{x}_1]}{N}
\biggr)g^{(2)*} \biggl(\frac{[N\mathbf
{y}]+\mathbf
{1}}{N},\frac{[Ns]\mathbf{1}-[N\mathbf{x}_2]}{N}
\biggr)
\nonumber
\\[-8pt]
\label{eqchaosproofbound}
\\[-8pt]
\nonumber
&&\qquad{}\times 1_{\{
[Ns]\mathbf
{1}>[N\mathbf{x}]\}\cap F(N)}
\\
\nonumber
&&\quad\le C g^{(1)*} (\mathbf{y},s\mathbf{1}-\mathbf{x}_1
)g^{(2)*} (\mathbf{y},s\mathbf{1}-\mathbf{x}_2
)1_{\{
s\mathbf
{1}>\mathbf{x}\}},
\end{eqnarray}
where we have used the following facts:
on the set $\{\mathbf{y}>\mathbf{0},[Ns]\mathbf{1}>[N\mathbf{x}]\}
$, we
have $([N\mathbf{y}]+1)/N>\mathbf{y}$, $([Ns]-[Nx_j])/N \ge\frac
{1}{2}(s-x_j)$ (see relation (40) in the proof of Theorem~6.5 of \cite{baitaqqu2013generalized}) and $g^{(j)*}$ decreases in its every variables,
as well as the fact that
$\{[Ns]\mathbf{1}>[N\mathbf{x}]\}\subset\{s\mathbf{1}>\mathbf{x}\}$.
Note that
%
\begin{eqnarray}
&&\int_0^t \mathrm{d}s \int
_{\mathbb{R}_+^r}\,\mathrm{d}\mathbf{y} g^{(1)*} (\mathbf {y},s\mathbf{1}-
\mathbf{x}_1 )g^{(2)*} (\mathbf {y},s\mathbf {1}-
\mathbf{x}_2 )1_{\{s\mathbf{1}>\mathbf{x}\}}
\nonumber
\\[-8pt]
\label{eqboundcontract}
\\[-8pt]
\nonumber
&&\quad=\int_0^t g^{(1)*}
\otimes_r g^{(2)*} (s\mathbf{1}-\mathbf{x} ) 1_{\{
s\mathbf{1}>\mathbf{x}\}}
\,\mathrm{d}s.
\end{eqnarray}
Since $g^{(1)*}$ and $g^{(2)*}$ are GHK(B)s, so by Lemma~\ref{LemgcontractgisGHK}, $g^{(1)*}\otimes_r g^{(2)*}$ is a GHK. This has two
consequences. First, by Theorem~3.5 and Remark~3.6 of \cite{baitaqqu2013generalized}, the integral in $\mathrm{d}s\,\mathrm{d}\mathbf{y}$ on the
left-hand side of (\ref{eqboundcontract}) is finite for a.e.
$\mathbf{x}\in\mathbb{R}^{k_1+k_2-2r}$. One can then apply the
dominated convergence theorem to conclude that
%
\begin{equation}
\label{eqgoalchaosproof} \tilde{h}_{N,t,r}(\mathbf{x})\rightarrow
h_{t,r}(\mathbf{x})\qquad\mbox{for a.e. }\mathbf{x}\in\mathbb{R}^{k_1+k_2-2r}.
\end{equation}
But to obtain (\ref{eqgoalcontractL^2}), we need $L^2$ convergence
for the integral in $\mathrm{d}\mathbf{x}$. For this, we use the bound~(\ref
{eqchaosproofbound}):
\[
\bigl\llvert \tilde{h}_{N,t,r}(\mathbf{x})\bigr\rrvert \le
h_{t,r}^*(\mathbf{x}):=C\int_0^t
g^{(1)*}\otimes_r g^{(2)*} (s\mathbf{1}-\mathbf{x} )
1_{\{
s\mathbf{1}>\mathbf{x}\}}\,\mathrm{d}s.
\]
The second consequence of the fact that $g^{(1)*}\otimes_r g^{(2)*}$ is
a GHK stems from Remark~\ref{RemGHKhtL^2}, which entails that
$h_{t,r}^*\in L^2(\mathbb{R}^{k_1+k_2-2r})$, and hence (\ref{eqgoalcontractL^2}) follows from (\ref{eqgoalchaosproof}) and the
dominated convergence theorem. This concludes the proof of Lemma~\ref
{Lemchaosproductnclt}.
\end{pf}

We now decompose
the product $X(n)$ in (\ref{eqchaosproductproc}) in off-diagonal
forms (\ref{eqXnnonsym}) as follows: let $\mathbf
{u}=(u_1,\ldots
,u_r)\in\mathbb{Z}_+^{r}$, $\mathbf{i}_1=(i_1,\ldots,i_{k_1-r})$ and
$\mathbf{i}_2=(i_{k_1-r+1},\ldots,i_{k_1+k_2-2r})$, and $\mathbf
{i}=(\mathbf{i}_1,\mathbf{i}_2)\in\mathbb{Z}_+^{k_1+k_2-2r}$, then
\begin{eqnarray*}
X(n)&=& Y'_1(n)Y'_2(n)\\
&=& \sum
_{r=0}^{k_1\wedge k_2}\! r! \pmatrix{k_1\cr r}
\pmatrix{k_2\cr r}\\
&&{}\times\sum_{(\mathbf{u},\mathbf{i})\in\mathbb
{Z}_+^{k_1+k_2-r}}'\! a^{(1)}(
\mathbf{u},\mathbf{i}_1)a^{(2)}(\mathbf {u},
\mathbf{i}_2) \varepsilon_{n-u_1}^2\cdots
\varepsilon_{n-u_r}^2 \varepsilon _{n-i_1}\cdots
\varepsilon_{n-i_{k_1+k_2-2r}}\!,
\end{eqnarray*}
where we have used the symmetry of $a^{(j)}$'s, while the combinatorial
coefficient
\[
c(r,k_1,k_2):=r! \pmatrix{k_1\cr r}
\pmatrix{k_2\cr r}
\]
is obtained as the number of ways to pair $r$ variables of $a^{(1)}$ to
$r$ variables of $a^{(2)}$.
We write
\[
\varepsilon_{n-i}^2=1+\bigl(\varepsilon_{n-i}^2-1
\bigr)=: A_0(\varepsilon _{n-i})+A_2(
\varepsilon_{n-i}),
\]
where $A_0(\varepsilon)=1$ and $A_2(\varepsilon)=\varepsilon^2-1$. These are
Appell polynomials which will be introduced in more details in Section~\ref{secproofvolt}. Set $J_r=\{0,2\}\times\cdots\times\{0,2\}$. Then
\begin{eqnarray*}
Y'_1(n)Y'_2(n)&= &\sum
_{r=0}^{k_1\wedge k_2} c(r,k_1,k_2)
\sum_{(\mathbf
{u},\mathbf{i})\in\mathbb{Z}_+^{k_1+k_2-r}}'\sum
_{\mathbf{j}\in J_r
}a^{(1)}(\mathbf{u},\mathbf{i}_1)a^{(2)}(
\mathbf{u},\mathbf{i}_2)
\\
&&{}\times A_{j_1}(\varepsilon_{n-u_1})\cdots A_{j_r}(
\varepsilon_{n-u_r}) \varepsilon _{n-i_1}\cdots\varepsilon_{n-i_{k_1+k_2-2r}}\!.
\end{eqnarray*}
The random variables in each summand are independent because the sum
does not include diagonals. Observe that it is only when $k_1=k_2$,
that the mean
\[
\mathbb{E}Y'_1(n)Y'_2(n)=k_1!
\sum_{\mathbf{u}\in\mathbb{Z}_+^{k_1}}' a^{(1)}(
\mathbf{u})a^{(2)}(\mathbf{u})
\]
may possibly be nonzero (this is the case when $r=k_1=k_2$). Hence, one
can use the $k$ defined in Theorem~\ref{Thmgeneralprodchaos} to
write that
%
\begin{eqnarray}
X(n)-\mathbb{E}X(n) &=&\sum_{r=0}^{k}\sum
_{\mathbf{j}\in J_r } \sum_{(\mathbf
{u},\mathbf{i})\in\mathbb{Z}_+^{k_1+k_2-r}}'
c(r,k_1,k_2) a^{(1)}(\mathbf{u},
\mathbf{i}_1)a^{(2)}(\mathbf{u},\mathbf {i}_2)
\nonumber
\\[-8pt]
\label{eqoffdiagdecompchaos}
\\[-8pt]
\nonumber
&&\hspace*{80pt}{}\times A_{j_1}(\varepsilon_{n-u_1})\cdots A_{j_r}(
\varepsilon_{n-u_r}) \varepsilon_{n-i_1}\cdots\varepsilon_{n-i_{k_1+k_2-2r}}.\quad
\end{eqnarray}
A basic term of the preceding decomposition of $X(n)-\mathbb{E}X(n)$ is
\begin{eqnarray*}
X_{\mathbf{j}}^r(n)&:=&\sum_{(\mathbf{u},\mathbf{i})\in\mathbb
{Z}_+^{k_1+k_2-r}}'c(r,k_1,k_2)
a^{(1)}(\mathbf{u},\mathbf {i}_1)a^{(2)}(
\mathbf{u},\mathbf{i}_2)
\\
&&\hspace*{38pt}\quad{}\times A_{j_1}(\varepsilon _{n-u_1})\cdots A_{j_r}(
\varepsilon_{n-u_r}) \varepsilon_{n-i_1}\cdots \varepsilon
_{n-i_{k_1+k_2-2r}}.
\end{eqnarray*}
Note\vspace*{1pt} that $0\le r\le k_1\wedge k_2$ if $k_1\neq k_2$, and $0\le r\le
k_1-1$ if $k_1=k_2$, which implies $k_1+k_2-2r\ge1$ so that there is
at least one $i$ variable. Due to the symmetry of $a^{(j)}$'s, we can
suppose without loss of generality that $j_1=\cdots= j_s=0$ and
$j_{s+1}=\cdots=j_r=2$, $0\le s\le r$. One can hence rewrite the basic
term as
\begin{eqnarray*}
X_{\mathbf{j}}^r(n)&=&\sum_{(\mathbf{u},\mathbf{i})\in\mathbb
{Z}_+^{k_1+k_2-r}}'c(r,k_1,k_2)
a^{(1)}(\mathbf{u},\mathbf {i}_1,\mathbf
{i}_2)a^{(2)}(\mathbf{u},\mathbf{i}_1,
\mathbf{i}_3)
\\
&&\hspace*{49pt}{}\times A_{2}(\varepsilon_{n-i_1})\cdots A_{2}(
\varepsilon_{n-i_{r-s}}) \varepsilon_{n-i_{r-s+1}}\cdots\varepsilon_{n-i_{k_1+k_2-r-s}},
\end{eqnarray*}
where
\begin{eqnarray*}
\mathbf{u} &=& (u_1,\ldots,u_s),\qquad \mathbf{i}_1=(i_1,
\ldots,i_{r-s}),
\\
\mathbf {i}_2 &=& (i_{r-s+1},\ldots ,i_{k_1-s}),\qquad
\mathbf{i}_3=(i_{k_1-s+1},\ldots,i_{k_1+k_2-r-s})
\end{eqnarray*}
and $\mathbf{i}=(\mathbf{i}_1,\mathbf{i}_2,\mathbf{i}_3)$. Setting
%
\begin{equation}
\label{eqchaosa} a'(\mathbf{i})=\sum_{\mathbf{u}\in K(\mathbf{i})}c(r,k_1,k_2)
a^{(1)}(\mathbf{u},\mathbf{i}_1,\mathbf{i}_2)a^{(2)}(
\mathbf {u},\mathbf {i}_1,\mathbf{i}_3),
\end{equation}
with
\[
K(\mathbf{i})=\{\mathbf{u}>\mathbf{0} \dvt u_p\neq u_q
\mbox{ if } p\neq q; \mbox{ and } u_p\neq i_q \mbox{
even if } p=q\},
\]
we get
%
\begin{equation}
\label{eqbasictermchaosproduct} X_{\mathbf{j}}^r(n)=\sum
_{\mathbf{i}>\mathbf{0}}'a'(\mathbf{i})
A_{2}(\varepsilon_{n-i_1})\cdots A_{2}(
\varepsilon_{n-i_{r-s}}) \varepsilon _{n-i_{r-s+1}}\cdots\varepsilon_{n-i_{k_1+k_2-r-s}}.
\end{equation}

We list here some useful elementary inequalities which will be used
many times in the sequel:

\begin{Lem} Let $A>0$, $B> 0$.  If $\gamma<-1$, then
%
\begin{equation}
\label{eqboundA+i} \sum_{i=1}^\infty(A+i)^{\gamma}
\le CA^{\gamma+1}.
\end{equation}
If $\gamma<0$, $\beta<-1$, then
%
\begin{equation}
\label{eqboundA+iismall} \sum_{i=1}^\infty(A+i)^{\gamma}i^{\beta}
\le CA^{\gamma}.
\end{equation}
If $\gamma<-1/2$, $-1<\beta<-1/2$, then
%
\begin{equation}
\label{eqboundA+iilarge} \sum_{i=1}^\infty(A+i)^{\gamma}i^\beta
\le CA^{\gamma+\beta+1}.
\end{equation}
If $\gamma<-1/2$, $\beta<-1/2$, then
%
\begin{equation}
\label{eqboundA+iB+ilarge} \sum_{i=1}^\infty(A+i)^{\gamma}(B+i)^{\beta}
\le C A^{\gamma
+1/2}B^{\beta+1/2}.
\end{equation}
\end{Lem}

\begin{pf}
To obtain inequality (\ref{eqboundA+i}), we have
\begin{eqnarray*}
\sum_{i=1}^\infty(A+i)^\gamma &=& \sum
_{i=1}^\infty\int_{i-1}^i
(A+i)^\gamma \,\mathrm{d}x \le\sum_{i=1}^\infty
\int_{i-1}^i (A+x)^\gamma \,\mathrm{d}x\\
&=& \int
_0^\infty(A+x)^\gamma \,\mathrm{d}x =-(
\gamma+1)^{-1}A^{\gamma+1}.
\end{eqnarray*}

For (\ref{eqboundA+iismall}), note that $(A+i)^\gamma\le
A^\gamma$
and $\sum_{i=1}^\infty i^{\beta}<\infty$.

For inequality (\ref{eqboundA+iilarge}),
we have
\[
\sum_{i=1}^\infty(A+i)^{\gamma}
i^{\beta} =A^{\gamma+\beta+1}\sum_{i=1}^\infty
\int_{i-1}^i (1+i/A)^{\gamma}
(i/A)^\beta \,\mathrm{d}(x/A) \le A^{\gamma+\beta+1} \int_0^\infty(1+y)^\gamma
y^{\beta}\,\mathrm{d}y,
\]
where the integral is finite since $\beta>-1$ and $\gamma+\beta<-1$.

The last one (\ref{eqboundA+iB+ilarge}) is obtained by applying
Cauchy--Schwarz and (\ref{eqboundA+i}) as follows:
\begin{eqnarray*}
\sum_{i=1}^\infty(A+i)^\gamma(B+i)^{\beta}
\le \Biggl[\sum_{i=1}^\infty
(A+i)^{2\gamma} \Biggr]^{1/2} \Biggl[ \sum
_{i=1}^\infty(B+i)^{2\beta
}
\Biggr]^{1/2} \le C A^{\gamma+1/2} B^{\beta+1/2}.
\end{eqnarray*}
\upqed\end{pf}
%
\begin{Rem}
The inequalities (\ref{eqboundA+i}), (\ref{eqboundA+iilarge})
and (\ref{eqboundA+iB+ilarge}) all raise the total power
exponent by $1$, while inequality (\ref{eqboundA+iismall}) kills
one of the exponents. These observations are useful in the proof below
and also in Section~\ref{secproofvolt}.
\end{Rem}

We now state the proof of Theorem~\ref{Thmgeneralprodchaos}.
\begin{pf*}{Proof of case~1 of Theorem~\ref
{Thmgeneralprodchaos}}
We want to apply Proposition~\ref{ThmcltXn}. The condition
$\mathbb{E}
|\varepsilon_i|^{4+\delta}<\infty$ guarantees that $\mathbb
{E}|A_2(\varepsilon
)|^{2+\delta'}<\infty$ in (\ref{eqbasictermchaosproduct}) holds for
some $\delta'>0$ and so the tightness in $D[0,1]$ holds.

We only need to show that $H^*<1/2$ in Lemma~\ref{Lembound} for each
of the basic terms $X_{\mathbf{j}}^r(n)$ in (\ref{eqbasictermchaosproduct}).

Suppose without loss of generality that $k_1\le k_2$.
Using the fact $|a^{(j)}(\mathbf{i})|\le C\|\mathbf{i}\|^{\alpha_j}$
(recall that $\|\cdot\|$ is the $L^1$-norm), one can bound $a'(\mathbf
{i})$ in (\ref{eqchaosa}). One has to distinguish two cases. In the
first case, where $s<k_1$, one gets
\begin{eqnarray*}
\bigl\llvert a'(\mathbf{i})\bigr\rrvert &\le & C \sum
_{\mathbf{u}\in\mathbb{Z}_+^s} \bigl\llVert (\mathbf{u},\mathbf {i}_1,
\mathbf{i}_2)\bigr\rrVert ^{\alpha_1} \bigl\llVert (\mathbf{u},
\mathbf{i}_1,\mathbf {i}_3)\bigr\rrVert ^{\alpha_2}
\\
&\le&  C\sum_{\mathbf{u}\in\mathbb{Z}_+^s} \bigl(u_1+
\cdots+u_s+\llVert \mathbf {i}_1\rrVert +\llVert
\mathbf{i}_2\rrVert \bigr)^{\alpha_1}\bigl(u_1+
\cdots+u_s+\llVert \mathbf{i}_1\rrVert +\llVert \mathbf
{i}_3\rrVert \bigr)^{\alpha_2}
\\
&\le &  C \bigl(\llVert \mathbf{i}_1\rrVert +\llVert
\mathbf{i}_2\rrVert \bigr)^{\alpha_1+s/2}\bigl(\llVert \mathbf
{i}_1\rrVert +\llVert \mathbf{i}_3\rrVert
\bigr)^{\alpha_2+s/2},
\end{eqnarray*}
after applying (\ref{eqboundA+iB+ilarge}) to each of the $s$
components of $\mathbf{u}$ iteratively (note: $\mathbf{i}_1$ may not be
present).
In the second case, where $s=r=k_1$, one gets
\begin{eqnarray*}
\bigl\llvert a'(\mathbf{i})\bigr\rrvert &\le  & C \sum
_{\mathbf{u}\in\mathbb{Z}_+^s} \llVert \mathbf{u}\rrVert ^{\alpha
_1} \bigl\llVert
(\mathbf{u},\mathbf{i}_3)\bigr\rrVert ^{\alpha_2}
\\
&\le &  C\sum_{\mathbf{u}\in\mathbb{Z}_+^s} (u_1+
\cdots+u_s)^{\alpha
_1}\bigl(u_1+\cdots+u_s+
\llVert \mathbf{i}_3\rrVert \bigr)^{\alpha_2}
\\
&\le & C \llVert \mathbf{i}\rrVert ^{\alpha_1+\alpha_2+s},
\end{eqnarray*}
after applying (\ref{eqboundA+iB+ilarge}) $s-1$ times, and then
(\ref{eqboundA+iilarge}) to the last component of $\mathbf{u}$. In
either case, the total power exponent is raised by $s$.

According to (\ref{eqH^*}), this yields
%
\begin{eqnarray}
H^*&=& \alpha_1+\alpha_2+s+(r-s+k_1-r+k_2-r)/2+1
\nonumber
\\
\label{eqH1+H2+s-r2-1}
&=& H_1+H_2+(s-r)/2-1
\\
&\le &  H_1+H_2-1<1/2,
\nonumber
\end{eqnarray}
where the last strict inequality is due to the assumption $H_1+H_2<3/2$
of case~1.
\end{pf*}
\begin{pf*}{Proof of case~2 of Theorem~\ref
{Thmgeneralprodchaos}}
We now suppose that $H_1+H_2>3/2$.
As was shown in case~1 above, the off-diagonal chaos
coefficient $a'(\cdot)$ in (\ref{eqchaosa}) leads to
\[
H^*=H_1+H_2+(s-r)/2-1.
\]

When $s=r$, we have only factors $A_0(\varepsilon)=1$ in (\ref{eqbasictermchaosproduct}). The chaos process $X_{r}^\mathbf{j}(n)$ is up to
some constant the process $X_r'(n)$ in Lemma~\ref{Lemchaosproductnclt}. Note that Lemma~\ref{Lemchaosproductnclt} concludes a joint
convergence for $X_r'(n)$ with different $r$'s. So adding up all the
terms corresponding to the case $r=s$ in (\ref{eqoffdiagdecompchaos}), which yields
\[
\sum_{r=0}^{k} \sum
_{(\mathbf{u},\mathbf{i})\in\mathbb
{Z}_+^{k_1+k_2-r}}' r! \pmatrix{k_1\cr r}
\pmatrix{k_2\cr r} a^{(1)}(\mathbf {u},\mathbf{i}_1)a^{(2)}(
\mathbf{u},\mathbf{i}_2) \varepsilon _{n-i_1}\cdots
\varepsilon_{n-i_{k_1+k_2-2r}},
\]
one obtains the noncentral limit claimed in the theorem with a Hurst
index $H=H_1+H_2-1>1/2$.

When $s<r$, the corresponding terms are negligible. Indeed,
\[
H^*=H_1+H_2+(s-r)/2-1\le H_1+H_2-1/2-1<
1/2.
\]
So by Lemma~\ref{Lembound}, the term $X_{r}^\mathbf{j}(n)$ has a
memory parameter $H\le1/2$ in the sense of Definition~\ref{DefHurstexpXn}. Hence,
\[
\lim_{N\rightarrow\infty}\mathbb{E} \Biggl[N^{-(H_1+H_2-1)}\sum
_{n=1}^{[Nt]}X_{r}^\mathbf{j}(n)
\Biggr]^2=0.
\]

We have now shown the convergence of finite-dimensional distributions.
Tightness in $D[0,1]$ is automatic since $H>1/2$ (see, e.g.,
Proposition~4.4.2 of \cite{giraitiskoulsurgailis2009large}).
\end{pf*}

\subsection{Proof of Theorem \texorpdfstring{\protect\ref{Thmgeneralprodvolt}}{3.6} where
diagonals are included}\label{secproofvolt}

We first recall from \cite{baitaqqu2014convergence} the
off-diagonal decomposition of a general $k$th order Volterra process
$X(n)$ in (\ref{eqXnintro}). The purpose is to decompose $X(n)$
into off-diagonal chaos terms as in (\ref{eqXnnonsym}).
To this end, it is convenient to use Appell polynomials. Suppose that
$\varepsilon$ is a random variable with finite $K$th moment. The Appell
polynomial with respect to the law of $\varepsilon$ is defined through the
following recursive relation:
\begin{eqnarray*}
\frac{\mathrm{d}}{\mathrm{d}x}A_p(x)= pA_{p-1},\qquad \mathbb{E}A_p(
\varepsilon)=0,\qquad A_0(x)=1,\qquad p=1,\ldots,K.
\end{eqnarray*}
We will use the following identity:
%
\begin{equation}
\label{eqappelldecomp} x^p = \sum_{j=0}^p
\pmatrix{p\cr j} \mu_{p-j}A_j(x),\qquad p=0,1,2,3,\ldots.
\end{equation}
For more details about Appell polynomials, see for example Chapter~3.3
of \cite{beran2013long}.

Let $\mathcal{P}_k$ be the collection of all the partitions of $\{
1,\ldots,k\}$. We further express each partition $\pi\in\mathcal{P}_k$
as $\pi=(P_1,\ldots,P_m)$ (so $m=|\pi|$), where the sets $P_t$'s are
\emph{ordered} according to their smallest element.
If we have a variable $\mathbf{i}\in\mathbb{Z}^k_+$, then $\mathbf
{i}_\pi$ denotes a new variable where its components are identified
according to $\pi$. For example, if $k=3$, $\pi=(\{1,2\},\{3\})$ and
$\mathbf{i}=(i_1,i_2,i_3)$, then $\mathbf{i}_\pi=(i_1,i_1,i_2)$. In
this case we write $\pi=(P_1,P_2)$ where $P_1=\{1,2\}$ and $P_2=\{3\}$.
If $a(\cdot)$ is a function on $\mathbb{Z}^k_+$, then
%
\begin{equation}
\label{eqapi} a_\pi(i_1,\ldots,i_m):=a(
\mathbf{i}_\pi),
\end{equation}
where $m=|\pi|$. In the preceding example, $a_{\pi}(\mathbf
{i})=a(i_1,i_2,i_2)$ with $m=2$. We define a summation operator $S'_T$
as follows:
for any $T\subset\{1,\ldots,|\pi|\}$, $S'_{T}(a_\pi)$ is obtained by
summing $a_\pi$ over its variables indicated by $T$ off-diagonally,
yielding a function with $|\pi|-|T|$ variables.
For instance, if $\pi=(\{1,5\},\{2\},\{3,4\})$, then $\mathbf{i}_\pi
=(i_1,i_2,i_3,i_3,i_1)$ and if $T=\{1,3\}$, then
\[
\bigl(S'_{T} a_\pi\bigr) (i)=\sum
_{0<i_1,i_3<\infty}' a(i_1,i,i_3,i_3,i_1),
\]
provided that it is well-defined. Note that in this off-diagonal sum,
we require also that neither $i_1$ nor $i_3$ equals to $i$.
If $T=\varnothing$, $S'_T$ is understood to be the identity operator.

Now, by collecting various diagonal cases and using (\ref{eqappelldecomp}), $X(n)$ in (\ref{eqXnintro}) can be decomposed as
%
\begin{eqnarray}
\label{eqoffdiagdecom} X(n)=\sum_{\pi\in\mathcal{P}_k} \sum
_{\mathbf{i}\in\mathbb
{Z}_+^m }' a_\pi(\mathbf{i})
\varepsilon_{n-i_1}^{p_1}\cdots\varepsilon _{n-i_m}^{p_m}=
\sum_{\pi\in\mathcal{P}_k}\sum_{\mathbf{j}\in
J(\pi)}
X_\pi^{\mathbf{j}}(n),
\end{eqnarray}
where
%
\begin{equation}
\label{eqbasicterm} X_\pi^{\mathbf{j}}(n)=\sum
_{\mathbf{i}\in\mathbb{Z}_+^m}' a_\pi (\mathbf {i})c(
\mathbf{p},\mathbf{j}) A_{j_1}(\varepsilon_{n-i_1})\cdots
A_{j_m}(\varepsilon_{n-i_m}),
\end{equation}
$A_j(\cdot)$ is the $j$th order Appell polynomial with respect to the
law of $\varepsilon_i$, $p_t=|P_t|$, $J(\pi)=\{0,\ldots,p_1\}\times
\cdots
\times\{0,\ldots,p_m\}$, and
%
\begin{equation}
\label{eqcpj} c(\mathbf{p},\mathbf{j})=\pmatrix{p_1\cr j_1}
\cdots \pmatrix{p_m\cr j_m}\mu _{p_1-j_1}\cdots
\mu_{p_m-j_m},\qquad\mu_j=\mathbb{E}\varepsilon_i^j.
\end{equation}
Note that since by assumption $\mu_1=0$, when $j_t=0$, it is only when
$p_t\ge2$ that it is possible to have a nonzero term.

In addition,
the expression for the centered $X(n)-\mathbb{E}X(n)$ is the sum in
(\ref
{eqoffdiagdecom}) with $J(\pi)$ replaced by $J^+(\pi):=J(\pi
)\setminus(0,\ldots,0)$, and
%
\begin{equation}
\label{eqmeanexpress} \mathbb{E}X(n)=\sum_{\pi\in\mathcal{P}_k}\sum
_{\mathbf{i}\in
\mathbb
{Z}^m_+}'a_{\pi}(\mathbf{i})
\mu_{p_1}\cdots\mu_{p_m}=\sum_{\pi
\in
\mathcal{P}_k^2}
\sum_{\mathbf{i}\in\mathbb{Z}^m_+}'a_{\pi
}(\mathbf
{i})\mu_{p_1}\cdots\mu_{p_m},
\end{equation}
where $\mathcal{P}_k^2$ denotes the collection of partitions of $\{
1,\ldots,k\}$ such that each set in the partition contains at least $2$
elements, namely, $p_t\ge2$ for all $t=1,\ldots,m$.

So from (\ref{eqoffdiagdecom}), (\ref{eqbasicterm}) and the
discussion above (\ref{eqmeanexpress}), the summands in the
off-diagonal decomposition of $X(n)-\mathbb{E}X(n)$ can be written as
%
\begin{equation}
\label{eqXpi^jbasicterm} X_\pi^{\mathbf{j}}(n)=\sum
_{\mathbf{i}\in\mathbb{Z}_+^{k'}}' c(\mathbf {p},\mathbf{j})
S_T'a_\pi(\mathbf{i}) A_{j_{t_1}}(
\varepsilon _{n-i_{t_1}})\cdots A_{j_{t_{k'}}}(\varepsilon_{n-i_{t_{k'}}}),
\end{equation}
where $T=\{t=1,\ldots,m \dvt j_t=0\}$, and $\{t_1,\ldots,t_{k'}\}= \{
1,\ldots,m\} \setminus T$ (thus $j_{t_1}\ge1,\ldots,j_{t_{k'}}\ge1$).
Note that $T\neq\{1,\ldots,m\}$ since $\mathbf{j}\in J^+(\pi)$. In
fact, $X_\pi^{\mathbf{j}}(n)$ is of the form (\ref{eqXnnonsym})
with $k=k'$ and $a(\cdot)=c(\mathbf{p},\mathbf{j}) S_T'a_\pi(\cdot)$.

We now state the proof of Theorem~\ref{Thmgeneralprodvolt} case by case.
Recall that $C>0$ denotes a constant whose value can change from line
to line.

\begin{pf*}{Proof of case~1}
In this case, $g^{(1)}(i)=C_1i^{\alpha_1}$, and
$g^{(2)}(i)=C_2i^{\alpha
_2}$, where $C_1$ and $C_2$ are two nonzero constants. The off-diagonal
decomposition (\ref{eqoffdiagdecom}) for the centered $X(n)$ is simply
%
\begin{equation}
\label{eqcase2simpleoffdiagdecomp} X(n)-\mathbb{E}X(n)= \sum_{0<i_1,i_2<\infty
}'a^{(1)}(i_1)a^{(2)}(i_2)
\varepsilon _{n-i_1}\varepsilon_{n-i_2}+\sum
_{0<i<\infty
}a^{(1)}(i)a^{(2)}(i)A_2(
\varepsilon_{n-i}),
\end{equation}
where $A_2(\varepsilon_{n-i})=\varepsilon_{n-i}^2-1$. Note that
\[
\bigl\llvert a^{(1)}(i_1)a^{(2)}(i_2)
\bigr\rrvert \le C i_1^{\alpha_1}i_2^{\alpha_2},
\]
so the coefficient of the first term in (\ref{eqcase2simpleoffdiagdecomp}) satisfies (\ref{eqabound}) with
\[
H^*=\alpha_1+\alpha_2+(1+1)/2+1=(H_1-3/2)+(H_2-3/2)+2<1/2
\]
by (\ref{eqHrange}), since $H_1+H_2<3/2$. For the second term in
(\ref
{eqabound}), one has
\[
\bigl\llvert a^{(1)}(i)a^{(2)}(i)\bigr\rrvert \le C
i^{\alpha_1+\alpha_2},
\]
which yields
%
\begin{equation}
\label{eqH*simplecase} H^*=\alpha_1+\alpha_2+1/2+1=(H_1-3/2)+(H_2-3/2)+3/2=H_1+H_2-3/2<1/2,
\end{equation}
since $H_1<1$ and $H_2<1$.
Hence,\vadjust{\goodbreak} Proposition~\ref{ThmcltXn} applies.
\end{pf*}

\begin{pf*}{Proof of case~2}
Now the first term of (\ref{eqcase2simpleoffdiagdecomp}) is
subject to Proposition~\ref{Thmncltchaossingle} with a Hurst index
$H=\alpha_1+\alpha_2+2=H_1+H_2-1>1/2$. One can see that for the second
term of (\ref{eqcase2simpleoffdiagdecomp}), relation (\ref{eqH*simplecase}) still holds. So by Lemma~\ref{Lembound}, the second term
of (\ref{eqcase2simpleoffdiagdecomp}) has a memory parameter
$H\le1/2$ in the sense of Definition~\ref{DefHurstexpXn}, and
hence with the normalization $N^{-H}$, the normalized partial sum of
the second term of (\ref{eqcase2simpleoffdiagdecomp}) converges
to $0$ in $D[0,1]$.
\end{pf*}

\begin{pf*}{Proof of case~3}
Recall from (\ref{eqXpi^jbasicterm}) that the summands in the
off-diagonal decomposition of $X(n)-\mathbb{E}X(n)$ are
\[
X_\pi^{\mathbf{j}}(n)=\sum_{\mathbf{i}\in\mathbb{Z}_+^{k'}} c(
\mathbf {p},\mathbf{j}) S_T'a_\pi(
\mathbf{i}) A_{j_{t_1}}(\varepsilon _{n-i_{t_1}})\cdots A_{j_{t_{k'}}}(
\varepsilon_{n-i_{t_{k'}}}).
\]
Consider first the following partition $\pi=(P_1,\ldots,P_m)$ of $\{
1,\ldots,k_1,k_1+1\}$, which we express as
\[
\pi=\bigl(P_1,\ldots,P_{m_1},\{k_1+1\}\bigr),
\]
with $m_1=m-1$, $\bigcup_{j=1}^{m_1} P_j=\{1,\ldots,k_1\}$, and $P_m=\{
k_1+1\}$. Let $T= \{1,\ldots,m_1\}$. Recall that to have nonzero
$c(\mathbf{p},\mathbf{j})$, one must require $|P_t|\ge2$ if $t\in T$,
and hence $2m_1\le k_1$. Set $\pi_1=\{P_1,\ldots,P_{m_1}\}$ and let
$\mathbf{u}\in\mathbb{Z}_+^{k_1}$. Then applying the off-diagonal
summation $S_T'$, we get
%
\begin{eqnarray}
\label{eqSTa} \bigl(S_T'a_{\pi}\bigr) (i)=
\sum_{u_p\neq u_q, u_p\neq i} a^{(1)}_{\pi
_1}(\mathbf
{u})a^{(2)}(i)= \biggl(\sum_{u_p\neq u_q}
a^{(1)}_{\pi_1}(\mathbf {u}) \biggr)a^{(2)}(i)-R(i),
\end{eqnarray}
where the difference $R(i)$ includes the terms where some $u_p=i$.
Since $|a^{(1)}(\mathbf{i})|\le C (i_1+\cdots+i_{k_1})^{\alpha_1}$
which implies $|a^{(1)}_{\pi}(\mathbf{u})|\le C (u_1+\cdots
+u_{m_1})^{\alpha_1}$. Suppose without loss of generality that
$u_{m_1}=i$, then by applying (\ref{eqboundA+i}),
\[
\bigl\llvert R(i)\bigr\rrvert \le C \sum_{0<u_1,\ldots,u_{m_1-1}<\infty}
(u_1+\cdots +u_{m_1-1}+i)^{\alpha_1} i^{\alpha_2} \le
C i^{\alpha_2+(\alpha_1+m_1-1)},
\]
where $\alpha_1+m_1-1<0$ because $\alpha_1<-k_1/2\le-m_1 \le-1$. It
follows that $|R(i)|\le Ci^{\alpha_2-\delta}$ for some $\delta>0$.
Since $k_2=1$,
the term $R(i)$ defines the linear process $\sum_{i>0} R(i) \varepsilon
_{n-i}$ but one with smaller memory parameter in the sense of
Definition~\ref{DefHurstexpXn}, than the linear process:
\[
\mu_{\pi_1} \Biggl(\sum_{\mathbf{u}>\mathbf{0}}'
a^{(1)}_{\pi
_1}(\mathbf {u}) \Biggr)\sum
_{i=1}^\infty a^{(2)}(i)
\varepsilon_{n-i},
\]
resulting from the first term in the right-hand side of (\ref{eqSTa}) (in this case $c(\mathbf{p},\mathbf{j})=\mu_{\pi_1}:=\mu
_{p_1}\cdots
\mu_{p_{m_1}}$).
Collecting all such $\pi_1\in\mathcal{C}_1^2$, one obtains $c_1\sum_{i=1}^\infty a^{(2)}(i) \varepsilon_{n-i}$ with $c_1$ as given in (\ref
{eqcsa}). Applying Proposition~\ref{Thmncltvoltsingle} with
$k=1$, we get the noncentral limit in case~3,
with a Hurst index
\[
H=\alpha_2+1/2+1=\alpha_2+3/2=H_2.
\]

We now show that in all the other cases, the memory parameter of
$X_{\pi
}^{\mathbf{j}}(n)$ is smaller than $H=\alpha_2+3/2$, which will
conclude the proof.
Observe first that
%
\begin{equation}
\label{eqaboundspecial} \bigl\llvert a(\mathbf{i})\bigr\rrvert \le C (i_1+
\cdots+i_{k_1})^{\alpha_1} i_{k_1+1}^{\alpha_2}.
\end{equation}
Let $\pi=\{P_1,\ldots,P_m\}$ is a partition of $\{1,\ldots,k_{1}+1\}$,
and $T=\{t_1,\ldots,t_l\}$, $l\le m-1$. To bound $|(S_T'a)(\mathbf{i})|$,
one can assume without loss of generality that either
\begin{longlist}[(a)]
\item[(a)] $P_j\cap\{k_1+1 \}=\varnothing$ for $1\le j\le m-1$, $P_m=\{
k_1+1\}$, $T\subset\{1,\ldots,m-1\}$, $\bigcup_{j=1}^{l} P_{t_j}\neq\{
1,\ldots,k_1\}$, or
\item[(b)] $P_m \cap\{k_1+1\}\neq\varnothing$, and $P_m\cap\{1,\ldots
,k_1\}\neq\varnothing$.
\end{longlist}
Observe that in the previous case we had $\bigcup_{j=1}^{l} P_{t_j}=\{
1,\ldots,k_1\}$ ($l=m_1=m-1$) and $P_m=\{k_1+1\}$.

In case (a), one has by (\ref{eqaboundspecial}) that
\[
\bigl\llvert a_\pi(\mathbf{i})\bigr\rrvert \le C(i_1+
\cdots+i_{m-1})^{\alpha
_1}i_{m}^{\alpha_2}.
\]
Since in case (a), $\bigcup_{j=1}^l P_{t_j}$ is a strict subset of $\{
1,\ldots,k_1\}$, we have $l< m-1$, and thus by applying (\ref{eqboundA+i}) iteratively, one has that
\begin{eqnarray*}
\bigl\llvert \bigl(S_T' a_\pi\bigr) (
\mathbf{i})\bigr\rrvert & \le& \sum_{\mathbf{u}>\mathbf{0}}
C(u_1+\cdots +u_l+i_{1}+
\cdots+i_{m-l-1})^{\alpha_1} i_{m-l}^{\alpha_2}
\\
&\le& C (i_{1}+\cdots+i_{m-l-1})^{\alpha_1+l}i_{m-l}^{\alpha_2},
\end{eqnarray*}
which results in $H^*$ in (\ref{eqH^*}) equal to
\begin{eqnarray*}
H^*&= &(\alpha_1+l+\alpha_2)+ (m-l)/2+1=
\alpha_1+\alpha_2+m/2+l/2+1
\\
&< & -k_1/2+\alpha_2 +(k_1+1)/2+1=
\alpha_2+3/2=H_2
\end{eqnarray*}
since $\alpha_1<-k_1/2$, and $m+l=2l+(m-l)\le k_1+1$ (recall that each
$|P_t|\ge2$ if $t\in T$).

In case (b), one can write without loss of generality that
\[
\bigl\llvert a_\pi(\mathbf{i})\bigr\rrvert \le C(i_1+
\cdots+i_m)^{\alpha_1}i_{1}^{\alpha_2}
\]
since $\pi$ contains $m$ partitions.
If for the above $a_\pi$, the summation $S_T'$ includes a sum over the
index $1$, that is, $1\in T$, then using (\ref{eqboundA+i}) and then
(\ref{eqboundA+iilarge}), one has
\begin{eqnarray*}
\bigl\llvert \bigl(S_T' a_\pi\bigr) (
\mathbf{i})\bigr\rrvert &\le & C \sum_{\mathbf{u}>\mathbf{0}}
(u_1+\cdots+u_l+i_1+\cdots+i_{m-l})^{\alpha_1}
u_{1}^{\alpha_2}
\\
&\le &  C \sum_{u_1=1}^\infty(u_1+
i_{1}+\cdots+i_{m-l})^{\alpha
_1+l-1}u_{1}^{\alpha_2}
\le C (i_{1}+\cdots+i_{m-l})^{\alpha_1+\alpha_2+l}.
\end{eqnarray*}
Relation (\ref{eqboundA+iilarge}) does apply because on one hand
$\alpha_2>-1$, and on the other hand, we have $\alpha_1+l-1<-1/2$ since
$\alpha_1<-k_1/2$ and $2(l-1)+1<k_1$ because of $|P_t|\ge2$ if $t\in
T$. This leads to $H^*$ in (\ref{eqH^*}) equal to
\begin{eqnarray*}
H^*&= &(\alpha_1+\alpha_2+l)+ (m-l)/2+1=
\alpha_1+\alpha _2+m/2+l/2+1<\alpha_2+3/2=H_2.
\end{eqnarray*}
If the summation $S_T'$ does not include the index $1$, that is, if
$1\notin T$, one has
\begin{eqnarray*}
\bigl\llvert \bigl(S_T' a_\pi\bigr) (
\mathbf{i})\bigr\rrvert &\le & C \sum_{\mathbf{u}> \mathbf{0}}
(i_1+\cdots+i_{m-l}+u_1+\cdots+u_{l})^{\alpha_1}
i_1^{\alpha_2}
\\
&\le  & C (i_{1}+\cdots+i_{m-l})^{\alpha_1+l}
i_1^{\alpha_2},
\end{eqnarray*}
by (\ref{eqboundA+i}),
which also yields $H^*<\alpha_2+3/2=H_2$.
\end{pf*}

\begin{pf*}{Proof of case~4}
Same as case~3.
\end{pf*}

\begin{pf*}{Proof of case~5}
We consider first in Part~1 all cases of $S_T'a_\pi$ in (\ref{eqXpi^jbasicterm}) which contribute to the limit, and in Part~2 negligible cases.

\emph{Part 1 of case~5}:
Suppose that $\pi$ can be split into $\pi_1$ and $\pi_2$ which satisfy
the following: the subpartition $\pi_1=\{P_1,\ldots,P_{m_1}\}$ is a
partition of $\{1,\ldots,k_1\}$, such that each $P_j$ satisfies
$|P_j|\le2$, and at least one $|P_j|=1$, $j=1,\ldots,m_1$.

Thus, suppose without loss of generality that $|P_1|=2,\ldots,|P_r|=2$,
$0\le r< m_1$, and $|P_{r+1}|=\cdots=|P_{m_1}|=1$. Require that the
subpartition $\pi_2$ belongs to $\mathcal{C}_2^2$, where $\mathcal
{C}_2^2$ is the collection of partitions of $\{k_1+1,\ldots,k_1+\cdots
+k_2\}$ such that each set in $\pi_2$ contains at least 2 elements.
$\mathcal{C}_2^2$~is nonempty because $k_2\ge2$. Let
\[
T=\{1,\ldots,r,m_1+1,\ldots,m_1+m_2\}.
\]
Setting $\mathbf{i}=(i_1,\ldots,i_{m_1-r})$, $\mathbf{u}=(u_1,\ldots
,u_r)\in\mathbb{Z}_+^r$ and $\mathbf{v}=(v_1,\ldots,v_{m_2})\in
\mathbb
{Z}_+^{m_2}$, one can write
%
\begin{eqnarray}\label{eqoridiagonalST}
\hspace*{-5pt}\bigl(S'_T a_\pi \bigr) (\mathbf{i})&=&
\mathop{\sum_{u_p\neq u_q u_p\neq i_q,u_p\neq v_q,}}_{v_p\neq v_q, v_p\neq i_q,\mathbf{u},\mathbf
{v}>\mathbf{0}}
a^{(1)}(u_1,u_1,\ldots,u_r,u_r,i_{1},
\ldots ,i_{m_1-r})a^{(2)}_{\pi_2}(\mathbf{v})
\\
\label{eqoridiagonalSTsecond}
&=&\sum_{u_p\neq u_q ,u_p\neq v_q,v_p\neq v_q,\mathbf
{u},\mathbf{v}>\mathbf{0}} a^{(1)}(u_1,u_1,
\ldots ,u_r,u_r,i_{1},\ldots
,i_{m_1-r})a^{(2)}_{\pi_2}(\mathbf{v})-R_1(
\mathbf{i})\qquad
\\
&=&\sum_{u_p\neq u_q,\mathbf{u}>\mathbf
{0}}a^{(1)}(u_1,u_1,
\ldots,u_r,u_r,i_{1},\ldots,i_{m_1-r})
\nonumber
\\[-8pt]
\label{eqtworesiduals}
\\[-8pt]
\nonumber
&&{}\times\sum_{v_p\neq
v_q,\mathbf{v}>\mathbf{0}}a^{(2)}_{\pi_2}(
\mathbf{v})-R_1(\mathbf {i})-R_2(\mathbf{i})
\end{eqnarray}
for $i_p\neq i_q$.
Relation (\ref{eqtworesiduals}) has the preceding three parts. We
shall now apply Proposition~\ref{Thmncltvoltsingle} to the first
part. Summing over all possible values of $r$, one gets a NCLT with
Hurst index $H=\alpha_1+k_1/2+1$, where the limit is
\[
Z:=c_2\sum_{0\le r<k_1/2} d_{k,r}
Z_{k_1-2r},
\]
where the process $Z_{k_1-2r}$ is defined in (\ref{eqZt=gr}) with
$g_r=g_r^{(1)}$. Taking into account that in this setting, $c(\mathbf
{p},\mathbf{j})$ in (\ref{eqcpj}) and (\ref{eqXpi^jbasicterm}) is
\begin{eqnarray*}
&& \pmatrix{p_1\cr 0}\cdots\pmatrix{p_r\cr 0} \pmatrix{p_{r+1} \cr 1} \cdots \pmatrix{p_{m_1} \cr 1} \pmatrix{p_{m_1+1} \cr 0}
\cdots \pmatrix{p_{m_1+m_2} \cr 0} (\mu _2)^r
\mu_{p_{m_1+1}}\cdots\mu_{p_{m_1+m_2}}
\\
&&\quad
=:\mu_{\pi_2},
\end{eqnarray*}
since $\mu_2=1$, $p_1=\cdots=p_r=2$ and $p_{r+1}=\cdots= p_{m_1}=1$,
one gets the nonzero constant $c_2$ in~(\ref{eqcsa}). As in (\ref
{eqwienerstratonintegral}), we can express the limit $Z(t)$ as a
centered Wiener--Stratonovich integral.

We shall now show that $R_1$ and $R_2$ in (\ref{eqtworesiduals}) lead
only to terms with Hurst indices strictly less than $H=\alpha
_1+k_1/2+1$ in the sense of Definition~\ref{DefHurstexpXn}, so
they are negligible compared to the first term, and hence they do not
contribute to the limit.

$R_1$ in (\ref{eqoridiagonalSTsecond}) is obtained by taking the
difference between the sum in (\ref{eqoridiagonalST}) and the sum
in (\ref{eqoridiagonalSTsecond}). Thus $R_1$ is obtained by
identifying some of the $u$ and $v$ variables in the sum in (\ref
{eqoridiagonalST}) with $i$ variables. Using the fact
$a^{(j)}(\mathbf{i})\le C\|\mathbf{i}\|^{\alpha_j}$, one can see that
one of the terms (a coefficient on $\mathbb{Z}_+^{m_1-r}$) in $R_1$ is
bounded by
%
\begin{equation}
\label{eqlastproof1} \sum_{u_p\neq u_q, u_p\neq v_q, v_p\neq v_q,\mathbf{u},\mathbf
{v}>\mathbf{0}} C\bigl(\llVert \mathbf{u}
\rrVert +\llVert \mathbf{i}\rrVert \bigr)^{\alpha_1}\bigl(\llVert \mathbf {v}
\rrVert +\bigl\llVert \mathbf{i}'\bigr\rrVert \bigr)^{\alpha_2},
\end{equation}
where $\mathbf{u}=(u_1,\ldots,u_{r-s_1})$, $\mathbf{i}=(i_1,\ldots
,i_{m_1-r})$, $\mathbf{v}=(v_1,\ldots,v_{m_2-s_2})$, $\mathbf
{i}'=(i_1,\ldots,i_{t})$, where
\[
0\le s_1 \le r\wedge(m_1-r),\qquad  0\le t\le s_2
\le m_2 \wedge(m_1-r).
\]
If $t=0$, then $s_2=0$, and in addition, either $s_1>0$ or $s_2>0$.
Note that $\mathbf{i}'$ is a subvector of $\mathbf{i}$.

By (\ref{eqboundA+i}), the term (\ref{eqlastproof1}) is bounded by
\begin{eqnarray*}
\sum_{\mathbf{u},\mathbf{v}>\mathbf{0}}C\bigl(\llVert \mathbf{u}\rrVert +\llVert
\mathbf {i}\rrVert \bigr)^{\alpha_1}\bigl(\llVert \mathbf{v}\rrVert +\bigl
\llVert \mathbf{i}'\bigr\rrVert \bigr)^{\alpha_2}\le %
\cases{ C\llVert \mathbf{i}\rrVert ^ {\alpha_1+r-s_1} &\quad$\mbox{if }t=0$;\vspace*{3pt}
\cr
C\llVert \mathbf{i}\rrVert ^ {\alpha_1+r-s_1} \bigl\llVert
\mathbf{i}'\bigr\rrVert ^ {\alpha_2+m_2-s_2} &\quad$\mbox{if }t>0$.}
\end{eqnarray*}
When $t=s_2=0$, one must have $s_1>0$, and so the term yields
\[
H^*=\alpha_1+r-s_1+(m_1-r)/2+1=
\alpha_1+(r+m_1)/2+1-s_1<\alpha_1+k_1/2+1,
\]
because
\[
r+m_1=2r+(m_1-r)=k_1.
\]

When $s_2\ge t>0$, it yields an
\begin{eqnarray*}
H^*&=&\alpha_1+r-s_1+\alpha_2+m_2-s_2+(m_1-r)/2+1
\\
&=&\alpha_1+(m_1+r)/2+1+\alpha_2+m_2-s_1-s_2
\\
&\le &\alpha_1+k_1/2+1+\alpha_2+k_2/2
-s_1-s_2<\alpha_1+k_1/2+1,
\end{eqnarray*}
since $2m_2\le k_2$ due to $\pi_2 \in\mathcal{C}_2^2$, and where the
last inequality is due to the assumption $\alpha_2<-k_2/2$.

We now examine $R_2$ in (\ref{eqtworesiduals}), which is obtained by
identifying some of the $u$ variables to the $v$ variables in the first
sum in (\ref{eqoridiagonalSTsecond}). One term of $R_2$ can be
bounded by
\[
\sum_{u_p\neq u_q, v_p\neq v_q,\mathbf{u},\mathbf{v}>\mathbf{0}} C \bigl(\llVert \mathbf{u}\rrVert +
\llVert \mathbf{v}_1\rrVert +\llVert \mathbf{i}\rrVert
\bigr)^{\alpha_1} \bigl(\llVert \mathbf {v}_1\rrVert +\llVert
\mathbf{v}_2\rrVert \bigr)^{\alpha_2},
\]
where $\mathbf{u}=(u_1,\ldots,u_{r-s}), \mathbf{v}_1=(v_1,\ldots,v_s),
\mathbf{v}_2=(v_{s+1},\ldots,v_{m_2})$ and $\mathbf{i}=(i_1,\ldots
,i_{m_1-r})$, where $1\le s\le(r\wedge m_2)$. By using (\ref{eqboundA+i}), and then (\ref{eqboundA+iB+ilarge}) and (\ref{eqboundA+iismall}), this term is bounded by
\begin{eqnarray*}
&&\sum_{\mathbf{u}>\mathbf{0},\mathbf{v}_1>\mathbf{0},\mathbf
{v}_2>\mathbf{0}} C \bigl(\llVert \mathbf{u}\rrVert +
\llVert \mathbf{v}_1\rrVert +\llVert \mathbf{i}\rrVert
\bigr)^{\alpha_1} \bigl(\llVert \mathbf{v}_1\rrVert +\llVert
\mathbf{v}_2\rrVert \bigr)^{\alpha_2}
\\
&&\quad\le \sum_{\mathbf{v}_1>\mathbf{0}} C \bigl(\llVert
\mathbf{v}_1\rrVert +\llVert \mathbf {i}\rrVert \bigr)^{\alpha_1+r-s}
\llVert \mathbf{v}_1\rrVert ^{\alpha_2+m_2-s} \le C \llVert \mathbf{i}
\rrVert ^{\alpha_1+r-s+(s-1)/2},
\end{eqnarray*}
which yields an
\[
H^*=\alpha_1+r-s/2-1/2+(m_1-r)/2+1=\alpha
_1+(m_1+r)/2+1-s/2-1/2<\alpha _1+k_1/2+1.
\]

So neither $R_1$ nor $R_2$ contributes to the limit.

\noindent\emph{Part 2 of case~5}.
Suppose now that $\pi$ and $T$ are \emph{not} as in Part 1. To
determine these cases, note that one can always bound $|(S_T'a_\pi
)(\mathbf{i})|$ by
%
\begin{equation}
\label{eqproofcasegeneral} C \sum_{\mathbf{u}>\mathbf{0}}\bigl(\llVert
\mathbf{i}_1\rrVert +\llVert \mathbf{i}_2\rrVert +\llVert
\mathbf{u}_1\rrVert +\llVert \mathbf{u}_2\rrVert
\bigr)^{\alpha_1}\bigl(\llVert \mathbf{i}_1\rrVert +\llVert
\mathbf {i}_3\rrVert +\llVert \mathbf{u}_1\rrVert +
\llVert \mathbf{u}_3\rrVert \bigr)^{\alpha_2},
\end{equation}
where $\mathbf{i}_j\in\mathbb{Z}_+^{s_j}$, $\mathbf{u}_j\in\mathbb
{Z}_+^{t_j}$, $s_j,t_j\ge0$ and where $s_1+s_2+s_3>0$ (at least one
$i$ variable must remain), and
\[
s_1+s_2+t_1+2t_2\le
k_1, \qquad s_1+s_3+t_1+2t_3
\le k_2.
\]
Thus, the variables in $\mathbf{u}_2$ are at least paired within
$a^{(1)}$, and the variables in $\mathbf{u}_3$ are at least paired
within $a^{(2)}$.

We note that in Part 1, we had $s_1=s_3=t_1=0$, and $s_t+2t_2=k_1$.
Thus, to avoid the situation considered in Part 1, we require
%
\begin{equation}
\label{eqavoidrep}
\mbox{if }s_1=s_3=t_1=0,\qquad
\mbox{then }s_2+2t_2<k_1.
\end{equation}

As we have dealt with $R_1$ and $R_2$ before, by properly applying
(\ref
{eqboundA+i})--(\ref{eqboundA+iB+ilarge}), the bound in (\ref
{eqproofcasegeneral}) yields
\[
H^*<H_1=\alpha_1/2+k_1/2+1.
\]
To check this, we consider the following exhaustive cases:
\begin{longlist}[(a)]
\item[(a)] either $s_1> 0$, or $s_1=0$, $s_2>0$, $s_3>0$;
\item[(b)]$s_1=s_2=0$, $s_3>0$;
\item[(c)]$s_1=s_3=0$, $s_2>0$ but $s_2+2t_2<k_1$.
\end{longlist}
Note that in case (c), if $s_2+2t_2=k_1$ then $t_1=0$, which would
contradict (\ref{eqavoidrep}).

In case (a), for example, if $s_1>0$, by applying (\ref{eqboundA+i})
to the sum over $\mathbf{u}_2$ and $\mathbf{u}_3$, and then (\ref
{eqboundA+iB+ilarge}) on the sum over $\mathbf{u}_1$, we can
bound (\ref{eqproofcasegeneral}) by
\begin{eqnarray*}
C \bigl(\llVert \mathbf{i}_1\rrVert +\llVert \mathbf{i}_2
\rrVert \bigr)^{\alpha_1+t_1/2+t_2}\bigl(\llVert \mathbf {i}_1\rrVert +
\llVert \mathbf{i}_3\rrVert \bigr)^{\alpha_2+t_1/2+t_3}.
\end{eqnarray*}
This yields
%
\begin{eqnarray}
H^*&=& \alpha_1+\alpha_2+t_1+t_2+t_3+(s_1+s_2+s_3)/2+1\nonumber
\\
\label{eqH*typical}
&=& \alpha_1+(s_1+s_2+t_1+2t_2)/2+1
+\alpha_2+(s_3+t_1+2t_3)/2
\\
&\le & \alpha_1+k_1/2+1 +\alpha_2+k_2/2<
\alpha_1+k_1/2+1=H_1.
\nonumber
\end{eqnarray}

In case (b), (\ref{eqproofcasegeneral}) becomes $C \sum_{\mathbf
{u}>\mathbf{0}} (\|\mathbf{u}_1\|+\|\mathbf{u}_2\|)^{\alpha_1}(\|
\mathbf
{i}_3\|+\|\mathbf{u}_1\|+\|\mathbf{u}_3\|)^{\alpha_2}$ which we can
bound by
\begin{eqnarray*}
&&C \sum_{\mathbf{u}_1>\mathbf{0}}\llVert \mathbf{u}_1
\rrVert ^{\alpha_1+t_2}\bigl(\llVert \mathbf{i}_3\rrVert +\llVert
\mathbf{u}_1\rrVert \bigr)^{\alpha_2+t_3}
\\
&&\quad\le %
\cases{ \llVert \mathbf{i}_3\rrVert
^{\alpha_2+(t_1-1)_+/2+t_3}& \quad$\mbox{ if } \alpha _1+t_1/2+t_2<-1/2$;
\vspace*{3pt}
\cr
\llVert \mathbf{i}_3\rrVert ^{\alpha_1+\alpha_2+t_1+t_2+t_3} & \quad$
\mbox{ if } -1/2<\alpha_1+t_1/2+t_2<0$,}
\end{eqnarray*}
where we need to apply first (\ref{eqboundA+i}), then apply (\ref
{eqboundA+iB+ilarge}) if $t_1\ge2$, and finally apply either
(\ref{eqboundA+iismall}) for the first case or (\ref{eqboundA+iilarge}) for the second. Note that $\alpha_1+t_1/2+t_2>-1/2$ only
if $t_1/2+t_2=k_1/2$ since $-k_1/2-1/2<\alpha_1<-k_1/2$ and
$t_1+2t_2\le k_1$. So this yields either an
%
\begin{eqnarray}
H^*&=&\alpha_2+(t_1-1)_+/2+t_3+s_3/2+1
\nonumber
\\[-8pt]
\label{eqcontributewhenH1=H2}
\\[-8pt]
\nonumber
&=&\alpha_2+ (s_3+t_1+2t_3)/2
+1 +(t_1-1)_+/2-t_1/2\le\alpha_2+k_2/2+1=
H_2<H_1
\end{eqnarray}
or $H^*$ as in (\ref{eqH*typical}).

Similarly, in case (c), (\ref{eqproofcasegeneral}) is $\sum_{\mathbf
{u}>\mathbf{0}}C(\|\mathbf{i}_2\|+\|\mathbf{u}_1\|+\|\mathbf{u}_2\|
)^{\alpha_1}(\|\mathbf{u}_1\|+\|\mathbf{u}_3\|)^{\alpha_2}$, which can
be bounded by
\begin{eqnarray*}
&&C \sum_{\mathbf{u}_1>\mathbf{0}}\bigl(\llVert \mathbf{i}_2
\rrVert +\llVert \mathbf{u}_1\rrVert \bigr)^{\alpha_1+t_2}\llVert
\mathbf{u}_1\rrVert ^{\alpha_2+t_3}
\\
&&\quad\le %
\cases{ \llVert \mathbf{i}_2\rrVert
^{\alpha_1+(t_1-1)_+/2+t_2}& \quad$\mbox{ if } \alpha _2+t_1/2+t_3<-1/2$;
\vspace*{3pt}
\cr
\llVert \mathbf{i}_2\rrVert ^{\alpha_1+\alpha_2+t_1+t_2+t_3} &\quad $
\mbox{ if } {-}1/2<\alpha_1+t_1/2+t_2<0$.}
\end{eqnarray*}
So it yields either an
%
\begin{eqnarray}
H^*&=&\alpha_1+(t_1-1)_+/2+t_2+s_2/2+1
\nonumber
\\[-8pt]
\label{eqH^*last}
\\[-8pt]
\nonumber
&=&\alpha_1+ (s_2+t_1+2t_2)/2 +1
+(t_1-1)_+/2-t_1/2<\alpha _1+k_1/2+1=H_1,
\end{eqnarray}
or $H^*$ as in (\ref{eqH*typical}).
To get the strict inequality in (\ref{eqH^*last}), we use (\ref
{eqavoidrep}) when $t_1=0$, and use $(t_1-1)_+/2<t_1/2$ when $t_1>0$.
\end{pf*}
\begin{pf*}{Proof of case~6}
Same as case~5.
\end{pf*}

\begin{pf*}{Proof of case~7}
Since $H_1=H_2$, both factors $a^{(1)}$ and $a^{(2)}$ may contribute to
the limit. The proof is similar to case~5, while the other term in the limit arises by exchanging of the
role of $a^{(1)}$ and $a^{(2)}$ in the proof of case~5. Note that because $H_1=H_2$, the equality in ``$\le
$'' in (\ref{eqcontributewhenH1=H2}) is attained whenever $t_1=0$
and $s_3+2t_3=k_2$, a case which would then be included in the NCLT
part of the proof.
\end{pf*}

\subsection{Proof of Theorem \texorpdfstring{\protect\ref{Thmgeneralprodmixed}}{3.8} the mixed
case}\label{secproofmixed}
%
The proof is similar to case 5 of Theorem~\ref{Thmgeneralprodvolt}.
We thus only give a sketch.

First, following the same notation as Part 1 of case 5 of Theorem~\ref
{Thmgeneralprodvolt}, we look at the contributing case: the
partition $\pi$ can be split into $\pi_1$ and $\pi_2$, where since now
the factor $Y_1'$ in (\ref{eqX=productmixed}) excludes the diagonals,
the first partition $\pi_1$ is just $\{\{1\},\ldots,\{m_1\}\}$. This
means that the
component $\mathbf{u}$ in (\ref{eqoridiagonalST}) does not appear,
namely, $r=0$. Hence, instead of (\ref{eqtworesiduals}) one gets
%
\begin{equation}
\label{eqmixedcontribute} a^{(1)}(i_{1},\ldots,i_{m_1})\sum
_{v_p\neq v_q,\mathbf{v}>\mathbf
{0}}a^{(2)}_{\pi_2}(\mathbf{v})-R(
\mathbf{i}),
\end{equation}
where $i_p\neq i_q$ for $p\neq q$ and the residual term $R(\mathbf{i})$
is as $R_1(\mathbf{i})$ in (\ref{eqtworesiduals}) (there is no $R_2$
due to absence of $\mathbf{u}$).
The first term leads to the noncentral limit $c_2 I_{k_1}(h_{t,1})$
with Hurst index $H_1=\alpha_1+k_1/2+1$ claimed in Theorem~\ref
{Thmgeneralprodmixed} by Proposition~\ref{Thmncltchaossingle}.
Then treating $R(\mathbf{i})$ in the same way as $R_1(\mathbf{i})$ is
treated there, one can show that $R(\mathbf{i})$ leads to terms with
Hurst index strictly less than $H_1=\alpha_1+k_1/2+1$. Since $H_1$ is
used in the normalization, all these terms are negligible.

Next, one follows Part 2 of case 5 of the proof of Theorem~\ref
{Thmgeneralprodvolt} to show that all other cases of $\pi$ yield
terms with Hurst indices strictly less than $H_1=\alpha_1+k_1/2+1$. Due
to the off-diagonality of $Y_1'$, for the bound (\ref{eqproofcasegeneral}), we have the following additional restrictions involving the
dimensions of the vectors in (\ref{eqproofcasegeneral}): $t_2=0$
($\mathbf{u}_2$ does not appear), and thus
%
\begin{equation}
\label{eqrestradd} s_1+s_2+t_1=k_1.
\end{equation}
The argument in the proof of Theorem~\ref{Thmgeneralprodvolt} for
cases (a) and (c) continue to hold because the quantity $H^*$ continues
to be strictly less than $H_1$. The only case there involving
modification is case (b) where $s_1=s_2=0, s_3>0$, because the original
inequality (\ref{eqcontributewhenH1=H2}) allows $H^*=H_2$ which
can be greater than $H_1$. But now by (\ref{eqrestradd}) we have the
restriction $t_1=k_1\ge2$. So (\ref{eqcontributewhenH1=H2}) is
now changed to
\[
H^*=\alpha_2+ (s_3+t_1+2t_3)/2
+1 -1/2\le H_2-1/2 < 1/2<H_1.
\]
Since $H^*<H_1$, these terms are also negligible. Then the first term
of (\ref{eqmixedcontribute}) dominates and provides the limit
$c_2I_{k_1}(h_{t,1})$.
%

\section*{Acknowledgements}

We would like to thank an anonymous
referee for careful reading and constructive suggestions. This work was
partially supported by the NSF Grants DMS-10-07616 and DMS-13-09009 at
Boston University.



%



\printhistory
\end{document}